\documentclass{amsart}

\theoremstyle{plain}
\newtheorem{theorem}{Theorem}[section]

\newtheorem{corollary}{Corollary}[section]
\newtheorem{proposition}{Proposition}[section]
\numberwithin{equation}{section}
\newcommand{\qbin}[2]{\genfrac{[}{]}{0pt}{}{#1}{#2}}
\newcommand{\qbins}[2]{{\textstyle\genfrac{[}{]}{0pt}{}{#1}{#2}}}

\begin{document}

\title[Partial theta functions]
{Partial theta functions. I. Beyond the lost notebook}

\author[Ole Warnaar]{S. Ole Warnaar}

\address{Department of Mathematics and Statistics,
The University of Melbourne, Vic 3010, Australia}
\email{warnaar@ms.unimelb.edu.au}


\thanks{Work supported by the Australian Research Council}

\begin{abstract}
It is shown how many of the partial theta function identities
in Ramanujan's lost notebook
can be generalized to infinite families of such identities.
Key in our construction is the Bailey lemma and a new 
generalization of the Jacobi triple product identity.
By computing residues around the poles of our identities
we find a surprising connection between partial theta functions 
identities and Garret--Ismail--Stanton-type
extensions of multisum Rogers--Ramanujan identities.
\end{abstract}

\maketitle

\section{Introduction}\label{sec1}

G.~E.~Andrews' discovery in 1976 of Ramanujan's lost 
notebook~\cite{Ramanujan88} can probably be regarded as one of the most
exciting finds ever in mathematics.
The lost notebook, which was hidden in a box containing papers from the late
G.~N.~Watson's estate, is a handwritten manuscript of well over a hundred
pages of hardly decipherable but very beautiful identities.
The first formula given by Andrews in his \textit{An introduction to 
Ramanujan's ``Lost'' notebook} \cite{Andrews79} is the following $q$-series
transformation \cite[p. 37]{Ramanujan88}:
\begin{multline}\label{R0}
\sum_{n=0}^{\infty}\frac{q^n}{\displaystyle (1-a)
\prod_{j=1}^n (1-aq^j)(1-q^j/a)} \\
=\sum_{n=0}^{\infty}(-1)^n a^{3n} q^{n(3n+1)/2}(1-a^2 q^{2n+1}) 
+\frac{\displaystyle \sum_{n=0}^{\infty}(-1)^n
a^{2n+1}q^{n(n+1)/2}}{\displaystyle \prod_{j=1}^{\infty}(1-aq^{j-1})(1-q^j/a)}
\end{multline}
on which Andrews puts the adjective ``marvelous''.

Characteristic of the above identity is that it contains
a partial theta product $(1-a)\prod_{j=1}^n (1-aq^j)(1-q^j/a)$
and a partial theta sum $\sum_{n=0}^{\infty} (-a^2)^n q^{n(n+1)/2}$.
Here it should be noted that complete theta products and sums are 
connected by Jacobi's famous triple product identity \cite[Eq. (2.28)]{GR90}
\begin{equation}\label{Jtpi}
\sum_{n=-\infty}^{\infty}(-1)^n a^n q^{n(n-1)/2}=
\prod_{n=1}^{\infty}(1-aq^{n-1})(1-q^n/a)(1-q^n).
\end{equation}

There are many more identities for ``partial theta functions''
in the lost notebook. With the standard notation
$(a;q)_n=\prod_{j=0}^{n-1}(1-aq^j)$ 
and $(a_1,\dots,a_k;q)_n =(a_1;q)_n\cdots(a_k;q)_n$, we state
the following further examples \cite[p. 12; p. 4; p. 29]{Ramanujan88}
\begin{gather}\label{R1}
\sum_{n=0}^{\infty}\frac{(q^{n+1};q)_n q^n}{(a;q)_{n+1}(q/a;q)_n}
=\sum_{n=0}^{\infty}a^n q^{n(n+1)}+\sum_{n=0}^{\infty}
\frac{a^{3n+1}q^{n(3n+2)}(1-aq^{2n+1})}{(a,q/a;q)_{\infty}}\\[3mm]
\label{R2}
\sum_{n=0}^{\infty}\frac{(q;q^2)_n q^n}{(a;q)_{n+1}(q/a;q)_n}
=\sum_{n=0}^{\infty}a^n q^{n(n+1)/2}+\sum_{n=0}^{\infty}
\frac{(-1)^n a^{2n+1}q^{n(n+1)}}{(-q,a,q/a;q)_{\infty}}
\end{gather}
and
\begin{multline}\label{R3}
\sum_{n=0}^{\infty}\frac{(q;q^2)_n q^{2n}}{(a;q^2)_{n+1}(q^2/a;q)_n} \\
=\sum_{n=0}^{\infty}a^n q^{n(n+1)/2}+
\sum_{n=0}^{\infty}\frac{(-1)^n a^{3n+1}q^{n(3n+2)}(1+aq^{2n+1})}
{(-q;q)_{\infty}(a,q^2/a;q^2)_{\infty}}.
\end{multline}

As was typical of Ramanujan, there are no proofs of any of
the partial theta function formulae in the lost notebook, 
making it impossible to determine how Ramanujan actually
discovered them.
Proofs of \eqref{R0} and \eqref{R1}--\eqref{R3} (and of many more
of Ramanujan's partial theta function identities) were found by Andrews 
\cite[Eq. (1.1), (1.2), (3.16) and (3.14)]{Andrews81}.
Equation \eqref{R0} was also proved by Andrews in \cite{Andrews79} and by
Fine \cite[Eqs. (7.2) and (7.5)]{Fine88}.
Andrews' proofs are at times quite intricate and rely heavily
on standard and some not-so-standard identities for
basic hypergeometric series. This perhaps partially explains why
Ramanujan's partial theta function identities, though beautiful 
and deep, have remained rather isolated and have not become as widely 
appreciated and studied as, for example, Ramanujan's mock theta functions.

The aim of this paper is to show that Ramanujan's partial theta function
identities are just the tip of the iceberg and that there is actually
a lot of hidden structure to \eqref{R0} and \eqref{R1}--\eqref{R3}.
For example, the identity \eqref{R1} is the first in an infinite
series of identities; the next identity in this series
being closely related to another of Ramanujan's discoveries,
that of the Rogers--Ramanujan identities \cite{Andrews76}. 
Specifically, we claim
\begin{multline*}
\sum_{n=0}^{\infty}
\frac{(q^{n+1};q)_n q^n}{(a;q)_{n+1}(q/a;q)_n}
\sum_{r=0}^n q^{r(r+1)} \qbin{n}{r}_q  \\
=\sum_{n=0}^{\infty} a^n q^{2n(n+1)}
+\sum_{i=1}^4(-1)^{i+1} a^i q^{\binom{i}{2}}
\frac{(q^i,q^{5-i},q^5;q^5)_{\infty}}{(q,a,q/a;q)_{\infty}}
\sum_{n=0}^{\infty} a^{5n}q^{2n(5n+2i)},
\end{multline*}
where the $q$-binomial coefficients or Gaussian polynomials are defined by
\begin{equation*}
\qbin{n}{m}_q=\begin{cases}\displaystyle \frac{(q;q)_n}{(q;q)_m(q;q)_{n-m}}&
\text{$m\in\{0,1,\dots,n\}$} \\[3mm]
0 & \text{otherwise.}
\end{cases}
\end{equation*}
More generally we have the following theorem.
\begin{theorem}\label{T1}
For $k$ a positive integer, $\kappa=2k+1$ and
$N_j=n_j+n_{j+1}+\cdots+n_{k-1}$ there holds
\begin{multline*}
\sum_{n=0}^{\infty}
\frac{(q;q)_{2n} q^n}{(a;q)_{n+1}(q/a;q)_n}
\sum_{n_1,\dots,n_{k-1}=0}^{\infty}
\frac{q^{N_1^2+\cdots+N_{k-1}^2+N_1+\cdots+N_{k-1}}}
{(q;q)_{n-N_1}(q;q)_{n_1}\cdots(q;q)_{n_{k-1}}}\\
=\sum_{n=0}^{\infty} a^n q^{kn(n+1)}
+\sum_{i=1}^{\kappa-1} (-1)^{i+1} a^i q^{\binom{i}{2}}
\frac{(q^i,q^{\kappa-i},q^{\kappa};q^{\kappa})_{\infty}}
{(q,a,q/a;q)_{\infty}}
\sum_{n=0}^{\infty} a^{\kappa n}q^{kn(\kappa n+2i)}.
\end{multline*}
\end{theorem}
Here we define $(a;q)_n$ for all integers $n$ as
\begin{equation}\label{aq}
(a;q)_n=\frac{(a;q)_{\infty}}{(aq^n;q)_{\infty}}
\end{equation}
so that, in particular, $1/(q;q)_n=0$ for negative $n$.
Observe that for $i\in\{1,2\}$ and $k=1$ the triple product 
$(q^i,q^{\kappa-i},q^{\kappa};q^{\kappa})_{\infty}$ becomes
$(q,q^2,q^3;q^3)_{\infty}=(q;q)_{\infty}$.
Since also $(q;q)_{2n}/(q;q)_n=(q^{n+1};q)_n$ one indeed finds
\eqref{R1} as the $k=1$ case of Theorem~\ref{T1}.

Our next theorem embeds \eqref{R2} in an infinite family.
\begin{theorem}\label{T2}
For $k$ a positive integer, $\kappa=2k$ and $N_j=n_j+n_{j+1}+\cdots+n_{k-1}$
there holds
\begin{align*}
\sum_{n=0}^{\infty} &
\frac{(q;q)_{2n} q^n}{(a;q)_{n+1}(q/a;q)_n}
\sum_{n_1,\dots,n_{k-1}=0}^{\infty}
\frac{q^{N_1^2+\cdots+N_{k-1}^2+N_1+\cdots+N_{k-1}}}
{(q;q)_{n-N_1}(q;q)_{n_1}\cdots(q;q)_{n_{k-2}}(q^2;q^2)_{n_{k-1}}} \\
&=\sum_{n=0}^{\infty} a^n q^{(\kappa-1)\binom{n+1}{2}} \\ 
& \qquad +\sum_{i=1}^{\kappa-1} (-1)^{i+1} a^i q^{\binom{i}{2}}
\frac{(q^i,q^{\kappa-i},q^{\kappa};q^{\kappa})_{\infty}}
{(q,a,q/a;q)_{\infty}}
\sum_{n=0}^{\infty} (-1)^n a^{\kappa n}q^{(\kappa-1)(kn+i)n}.
\end{align*}
\end{theorem}
To extend \eqref{R3} we have to rewrite the term $(q;q^2)_n$ in the
summand on the left as $(q^2;q^2)_{2n}/(q^2;q^2)_n(-q;q)_{2n}$.
\begin{theorem}\label{T3}
For $k$ a positive integer, $\kappa=2k-1/2$ and
$N_j=n_j+n_{j+1}+\cdots+n_{k-1}$ there holds
\begin{align*}
\sum_{n=0}^{\infty} &
\frac{(q;q)_{2n} q^n}{(a;q)_{n+1}(q/a;q)_n}
\sum_{n_1,\dots,n_{k-1}=0}^{\infty}
\frac{q^{N_1^2+\cdots+N_{k-1}^2+N_1+\cdots+N_{k-1}}}
{(q;q)_{n-N_1}(q;q)_{n_1}\cdots(q;q)_{n_{k-1}}(-q^{1/2};q^{1/2})_{2n_{k-1}}}\\
&=\sum_{n=0}^{\infty} a^n q^{(\kappa-1)\binom{n+1}{2}} 
+\sum_{i=1}^{2k-1} (-1)^{i+1} a^i q^{\binom{i}{2}}
\frac{(q^i,q^{\kappa-i},q^{\kappa};q^{\kappa})_{\infty}}
{(q,a,q/a;q)_{\infty}} \\
& \qquad \qquad \qquad \times
\sum_{n=0}^{\infty} (-1)^n a^{2\kappa n}q^{2(\kappa-1)(\kappa n+i)n}
\bigl(1+a^{2\kappa-2i}q^{2(\kappa-1)(\kappa-i)(2n+1)}\bigr).
\end{align*}
\end{theorem}
Finally, the generalization of \eqref{R0} to an infinite series is 
more complicated, involving a quintuple instead of triple product.
\begin{theorem}\label{T4}
For $k$ a positive integer, $\kappa=3k-1$ and
$N_j=n_j+n_{j+1}+\cdots+n_{k-1}$ there holds
\begin{align*}
\sum_{n=0}^{\infty}&
\frac{(q;q)_{2n} q^n}{(a;q)_{n+1}(q/a;q)_n}
\sum_{n_1,\dots,n_{k-1}=0}^{\infty}
\frac{q^{N_1^2+\cdots+N_{k-1}^2+N_1+\cdots+N_{k-1}}}
{(q;q)_{n-N_1}(q;q)_{n_1}\cdots(q;q)_{n_{k-2}}(q;q)_{2n_{k-1}}} \\
&=\sum_{n=0}^{\infty} (-1)^n a^{3n} q^{(2\kappa-3)(3n+1)n/2}
\bigl(1-a^2 q^{(2\kappa-3)(2n+1)}\bigr) \\
& \quad +\sum_{i=1}^{\kappa-1}(-1)^{i+1} a^i q^{\binom{i}{2}}
\frac{(q^i,q^{2\kappa-i},q^{2\kappa};q^{2\kappa})_{\infty}
(q^{2\kappa-2i},q^{2\kappa+2i};q^{4\kappa})_{\infty}}{(q,a,q/a;q)_{\infty}} \\
&\qquad \qquad \qquad \times \Bigl[1-\sum_{n=1}^{\infty}
a^{2\kappa n-2i}q^{(2\kappa-3)(\kappa n-i)n}
\bigl(1-a^{2i}q^{2(2\kappa-3)in}\bigr)\Bigr].
\end{align*}
\end{theorem}
To see this generalizes \eqref{R0}, note that
$(q^i,q^{2\kappa-i},q^{2\kappa};q^{2\kappa})_{\infty}
(q^{2\kappa-2i},q^{2\kappa+2i};q^{4\kappa})_{\infty}$ for
$k=i=1$ becomes $(q,q^3,q^4;q^4)_{\infty}(q^2,q^6;q^8)_{\infty}=
(q,q^2,q^3,q^4;q^4)_{\infty}=(q;q)_{\infty}$.

Several further results similar to Theorems~\ref{T1}--\ref{T4},
but reducing for $k=1$ to partial theta function identities not in the
lost notebook, will also be proved in this paper.
Apart from some deep but known results from the theory
of $q$-series, the following generalized triple product identity will
be crucial in our derivation of partial theta function identities.
\begin{theorem}\label{thmpts}
There holds
\begin{equation}\label{pts}
1+\sum_{n=1}^{\infty}(-1)^n q^{n(n-1)/2}(a^n+b^n)
=(q,a,b;q)_{\infty}\sum_{n=0}^{\infty}
\frac{(ab/q;q)_{2n}q^n}{(q,a,b,ab;q)_n}.
\end{equation}
\end{theorem}
We believe this to be a very beautiful formula.
Note in particular that for $b=q/a$ one recovers the Jacobi triple
product identity \eqref{Jtpi}. Indeed, making this specialization
only the $n=0$ term contributes to the sum on the right, whereas the
left simplifies to $\sum_{n=-\infty}^{\infty} a^n q^{n(n-1)/2}$.
Another nice specialization occurs when $b=-a$. Substituting this
and replacing $a$ by $(aq)^{1/2}$ gives the transformation
\begin{equation*}
1+2\sum_{n=1}^{\infty}a^n q^{2n^2}=(q;q)_{\infty}(aq;q^2)_{\infty}
\sum_{n=0}^{\infty} \frac{(-a;q)_{2n}q^n}{(q,-aq;q)_n(aq;q^2)_n}.
\end{equation*}

The remainder of this paper is organized as follows.
The next section contains a proof of the key identity \eqref{pts}.
In section~\ref{sec3} it is shown how \eqref{pts} can be applied
to give a very general partial theta function identity.
As examples Ramanujan's identities \eqref{R1} and \eqref{R2}
are obtained. Section~\ref{sec4} is devoted to the
proof of Theorems~\ref{T1}--\ref{T4} and related identities.
In sections~\ref{sec5} and \ref{sec6} we exploit the fact that all of the
partial theta function identities exhibit poles at $a=q^N$.
Calculating the residues around these simple poles
yields new identities which turn out to be Rogers--Ramanujan identities
of the Garret--Ismail--Stanton type.
We conclude the paper in section~\ref{sec7}
with a brief discussion of the possibilities and
limitations of our approach to partial theta function identities. 
The proofs of several polynomial identities needed in the main text
can be found in an appendix.

\section{Proof of Theorem \ref{thmpts}}
As a first step we use \eqref{aq} and $(a;q)_{2n}=(a;q)_n(aq^n;q)_n$ to
put \eqref{pts} in the form
\begin{multline}\label{pts2}
\sum_{n=0}^{\infty}
\frac{(aq^n,bq^n;q)_{\infty}(abq^{n-1};q)_n q^n}
{(1-abq^{n-1})(q;q)_n} \\
=\frac{1}{(1-ab/q)(q;q)_{\infty}}\Bigl\{1+\sum_{n=1}^{\infty}(-1)^n
q^{\binom{n}{2}}(a^n+b^n)\Bigr\}.
\end{multline}
Expanding the left side using the $q$-binomial theorem
\cite[Eq. (3.3.6)]{Andrews76}
\begin{equation}\label{qbt}
\sum_{k=0}^n (-z)^k q^{\binom{k}{2}}\qbin{n}{k}_q=(z;q)_n
\end{equation}
and its limiting $q$-exponential sum \cite[Eq. (II.2)]{GR90}
\begin{equation}\label{qE}
\sum_{k=0}^{\infty}\frac{(-z)^k q^{\binom{k}{2}}}{(q;q)_k}
=(z;q)_{\infty}
\end{equation}
gives the quadruple sum
\begin{equation*}
\text{LHS}\eqref{pts2}=
\sum_{n,i,j,k,l=0}^{\infty}
\frac{(-1)^{i+j+k}a^{i+k+l}b^{j+k+l}q^{\binom{i}{2}+\binom{j}{2}+
\binom{k}{2}+n(i+j+k+l+1)-k-l}}
{(q;q)_i(q;q)_j(q;q)_k(q;q)_{n-k}}.
\end{equation*} 
Shifting $n\to n+k$, then summing over $n$ using \cite[Eq. (II.1)]{GR90}
\begin{equation}\label{qe}
\sum_{n=0}^{\infty} \frac{z^n}{(q;q)_n}=\frac{1}{(z;q)_{\infty}},
\end{equation}
and finally shifting $i\to i-k$ and $j\to j-k$ gives
\begin{equation*}
\text{LHS}\eqref{pts2}=
\sum_{\substack{i,j,k,l=0\\ k\leq \min\{i,j\}}}^{\infty}
\frac{(-1)^{i+j+k}a^{i+l}b^{j+l}q^{\binom{i}{2}+\binom{j}{2}+
\binom{k+1}{2}+l(k-1)}(q;q)_{i+j+l-k}}
{(q;q)_{\infty}(q;q)_{i-k}(q;q)_{j-k}(q;q)_k}.
\end{equation*} 
Here the condition on the sum over $k$ is added to avoid possible ambiguity
for $i+j+l-k<0$. Employing the standard $q$-hypergeometric notation \cite{GR90}
\begin{align*}
{_{r+1}\phi_r}\Bigl[\genfrac{}{}{0pt}{}
{a_1,\dots,a_{r+1}}{b_1,\dots,b_r};q,z\Bigr]&=
{_{r+1}\phi_r}(a_1,\dots,a_{r+1};b_1,\dots,b_r;q,z) \\
&=\sum_{n=0}^{\infty}
\frac{(a_1,\dots,a_{r+1};q)_n}{(q,b_1,\dots,b_r;q)_n}z^n
\end{align*}
we can carry out the sum over $k$
by the $q$-Chu--Vandermonde sum \cite[(II.6)]{GR90}
\begin{equation}\label{qCV}
{_2\phi_1}(a,q^{-n};c;q,q)=\frac{(c/a;q)_n}{(c;q)_n}a^n,
\end{equation}
with $a=q^{-i}$, $n=j$ and $c=q^{-i-j-l}$. This yields
\begin{equation*}
\text{LHS}\eqref{pts2}=\sum_{i,j,l=0}^{\infty}
\frac{(-1)^{i+j}a^{i+l}b^{j+l}q^{\binom{i}{2}+\binom{j}{2}-l}
(q;q)_{i+l}(q;q)_{j+l}}{(q;q)_{\infty}(q;q)_i(q;q)_j(q;q)_l}.
\end{equation*}
After the changes $i\to i-l$ and $j\to j-l$ the sum over
$l$ can again be performed by \eqref{qCV}, now with $a=q^{-i}$, $n=j$ 
and $c=0$. Hence
\begin{equation*}
\text{LHS}\eqref{pts2}=\frac{1}{(q;q)_{\infty}}\sum_{i,j=0}^{\infty}
(-1)^{i+j}a^i b^j q^{\binom{i}{2}+\binom{j}{2}-i j}.
\end{equation*}
Equating this with the right-hand side of \eqref{pts2} we are left to
show that
\begin{equation}\label{ijn}
(1-ab/q)
\sum_{i,j=0}^{\infty}
(-1)^{i+j}a^i b^j q^{\binom{i}{2}+\binom{j}{2}-i j}
=1+\sum_{n=1}^{\infty}(-1)^n q^{\binom{n}{2}}(a^n+b^n).
\end{equation}
Making the change $i\to i+j$ on the left and taking care to respect the
ranges of summation gives
\begin{align*}
\text{LHS}\eqref{ijn}&=
(1-ab/q)\sum_{i=-\infty}^{\infty}\sum_{j=\max\{0,-i\}}^{\infty} 
(-1)^i a^i (ab/q)^j q^{\binom{i}{2}}  \\
&=\sum_{i=-\infty}^{\infty} 
(-1)^i a^i (ab/q)^{\max\{0,-i\}}q^{\binom{i}{2}}  \\
&=1+\sum_{i=1}^{\infty} (-1)^i a^i q^{\binom{i}{2}}
+\sum_{i=-\infty}^{-1}(-1)^i b^{-i} q^{\binom{i+1}{2}}=\text{RHS}\eqref{ijn}.
\end{align*}

\section{A general partial theta function identity}\label{sec3}
In this section we will show how \eqref{pts} can be applied to
yield a very general identity for partial theta functions
from which Ramanujan's identities of the introduction easily follow.
First we replace $a\to aq^{r+1}$ and $b\to bq^r$ in \eqref{pts}.
Using \eqref{aq}, $(aq^k;q)_n=(a;q)_{n+k}/(a;q)_k$ and shifting 
$n\to n-r$ on the right this can be written as
\begin{multline*}
1+\sum_{n=1}^{\infty}(-1)^n q^{\binom{n}{2}}
\bigl((aq^{r+1})^n+(bq^r)^n\bigr) \\
=q^{-r}(q,a,b;q)_{\infty}\frac{1-abq^{2r}}{1-ab}
\sum_{n=r}^{\infty} \frac{(ab;q)_{2n}q^n}{(q;q)_{n-r}(abq;q)_{n+r}
(a;q)_{n+1}(b;q)_n}.
\end{multline*}
By the use of the triple product identity \eqref{Jtpi} it follows that
\begin{multline*}
(-1)^r q^{\binom{r+1}{2}} \Bigl\{1+\sum_{n=1}^{\infty}
(-1)^n q^{\binom{n}{2}}\bigl((aq^{r+1})^n+(bq^r)^n\bigr) \Bigr\} \\
=(q/b)^r (q,b,q/b;q)_{\infty}
+\sum_{n=1}^{\infty}(-1)^{n+r} \Bigl\{ a^n q^{\binom{n+r+1}{2}}-
(q/b)^n q^{\binom{n-r}{2}} \Bigr\}. \\
\end{multline*}
Consequently there holds
\begin{multline*}
(q/b)^r (q,b,q/b;q)_{\infty}
+\sum_{n=1}^{\infty}(-1)^{n+r} \Bigl\{a^n q^{\binom{n+r+1}{2}}
-(q/b)^n q^{\binom{n-r}{2}}\Bigr\} \\
=(-1)^r q^{\binom{r}{2}}(q,a,b;q)_{\infty}\frac{1-abq^{2r}}{1-ab}
\sum_{n=r}^{\infty} \frac{(ab;q)_{2n}q^n}{(q;q)_{n-r}(abq;q)_{n+r}
(a;q)_{n+1}(b;q)_n}.
\end{multline*}
Next we multiply both sides by $f_r$ and sum $r$ over
the nonnegative integers leading to
\begin{multline}\label{ab}
\sum_{n=0}^{\infty}\frac{(ab;q)_{2n}q^n}{(a;q)_{n+1}(b;q)_n}
\sum_{r=0}^n \frac{(-1)^r q^{\binom{r}{2}}f_r(1-abq^{2r})}
{(q;q)_{n-r}(ab;q)_{n+r+1}}
-\frac{(q/b;q)_{\infty}}{(a;q)_{\infty}}
\sum_{r=0}^{\infty}(q/b)^r f_r \\
=\frac{1}{(q,a,b;q)_{\infty}}\sum_{n=1}^{\infty}\Bigl\{
a^n \sum_{r=0}^{\infty}f_r+(q/b)^n \sum_{r=-\infty}^{-1}f_{-r-1}\Bigr\}
(-1)^{n+r} q^{\binom{n+r+1}{2}}.
\end{multline}
Of course one either has to impose conditions on $f_r$ to 
ensure convergence of all sums or one has to view this identity
as a formal power series in $q$. 
Nearly all our applications of this general result assume the relation
$b=q/a$.
\begin{proposition}\label{prop1}
As a formal power series in $a$ and $q$ there holds
\begin{multline}\label{eqprop}
\sum_{n=0}^{\infty}\frac{(q;q)_{2n}q^n}{(a;q)_{n+1}(q/a;q)_n}
\sum_{r=0}^n \frac{(-1)^r q^{\binom{r}{2}}f_r(1-q^{2r+1})}
{(q;q)_{n-r}(q;q)_{n+r+1}}
-\sum_{r=0}^{\infty}a^r f_r \\
=\frac{1}{(q,a,q/a;q)_{\infty}}\sum_{n=1}^{\infty}\Bigl\{
\sum_{r=0}^{\infty}f_r+\sum_{r=-\infty}^{-1}f_{-r-1}\Bigr\}
(-1)^{n+r} a^n q^{\binom{n+r+1}{2}}.
\end{multline}
\end{proposition}
All that remains to be done to turn this into a Ramanujan-type
partial theta function identity is to appropriately choose
$f_r$ such that the sum 
\begin{equation}\label{beta}
\beta_n=\sum_{r=0}^n \frac{(-1)^r q^{\binom{r}{2}}f_r(1-q^{2r+1})}
{(q;q)_{n-r}(q;q)_{n+r+1}}
\end{equation}
can be carried out explicitly.
The most general such $f_r$ appears to be
\begin{equation}\label{f}
f_r=\frac{(b,c;q)_r}{(q^2/b,q^2/c;q)_r}\Bigl(\frac{q^2}{bc}\Bigr)^r.
\end{equation}
Then, by the $a=q$ case of Rogers' $q$-Dougall sum \cite[Eq. (II.21)]{GR90}
\begin{equation}\label{qDougall}
{_6\phi_5}\Bigl[\genfrac{}{}{0pt}{}
{a,qa^{1/2},-qa^{1/2},b,c,q^{-n}}{a^{1/2},-a^{1/2},aq/b,aq/c,aq^{n+1}};q,
\frac{aq^{n+1}}{bc}\Bigr]=\frac{(aq,aq/bc;q)_n}{(aq/b,aq/c;q)_n},
\end{equation}
it follows that
\begin{equation}\label{betabc}
\beta_n=\frac{(q^2/bc;q)_n}{(q,q^2/b,q^2/c;q)_n}.
\end{equation}
Inserted in Proposition~\ref{prop1} this leads to
\begin{multline}\label{bc}
\sum_{n=0}^{\infty}
\frac{(q^{n+1},q^2/bc;q)_n q^n}{(a;q)_{n+1}(q/a,q^2/b,q^2/c;q)_n}
-\sum_{r=0}^{\infty}\frac{(b,c;q)_r}{(q^2/b,q^2/c;q)_r}
\Bigl(\frac{aq^2}{bc}\Bigr)^r \\
=\sum_{n=1}^{\infty}
\frac{(-1)^n a^n q^{\binom{n+1}{2}}}{(q,a,q/a;q)_{\infty}}
\sum_{r=-\infty}^{\infty} 
\frac{(-1)^r q^{\binom{r}{2}}(b,c;q)_r}{(q^2/b,q^2/c;q)_r}
\Bigl(\frac{q^{n+3}}{bc}\Bigr)^r, 
\end{multline}
where we have used the symmetry $f_{-r-1}=f_r$ as follows from
\begin{equation*}
(a;q)_{-n}=\frac{(-q/a)^n q^{\binom{n}{2}}}{(q/a;q)_n}.
\end{equation*}
The partial theta function identities \eqref{R1} and \eqref{R2}
of Ramanujan arise as special limiting cases of this identity.
First, when $b$ and $c$ tend to infinity the sum over $r$ on the 
right simplifies to
\begin{align*}
\sum_{r=-\infty}^{\infty} 
(-1)^r q^{3\binom{r}{2}}q^{(n+3)r}
&=(q^{n+3},q^{-n},q^3;q^3)_{\infty} \\
&=(q;q)_{\infty}(-1)^{\lfloor (n+2)/3\rfloor}q^{-(n+1)(n+2)/6}
\chi(n\not\equiv 0 \pmod{3}),
\end{align*}
where we have used the triple product identity \eqref{Jtpi}
and 
\begin{equation}\label{pred}
(q^{i+mn},q^{m-mn-i};q^m)_{\infty}=
(q^i,q^{m-i};q^m)_{\infty}(-1)^n q^{-ni-\binom{n}{2}m},
\end{equation}
and where $\chi(\text{true})=1$, $\chi(\text{false})=0$ and
$\lfloor x \rfloor$ is the integer part of $x$.
Also letting $b,c\to\infty$ in the other terms in \eqref{bc}
we obtain Ramanujan's \eqref{R1}. This solves a problem of Andrews 
who remarked in \cite{Andrews81}: 
``The primary reason that our proof is so complicated is that we have
been unable to prove any generalization of (1.2)$_R$''.
Here (1.2)$_R$ is our \eqref{R1} and the type of generalization
Andrews alludes to is not a generalization like Theorem~\ref{T1}
but a generalization involving additional free parameters.

In much the same way as we obtained \eqref{R1} one finds \eqref{R2} 
after taking $b=-q$ and $c\to\infty$ in \eqref{bc}.

The $q$-Dougall sum \eqref{qDougall} can also be used to
derive quadratic and cubic analogues of \eqref{bc}.
Specifically, from \eqref{qDougall} it follows that
\begin{equation*}
f_{2r}=\frac{(q,b;q^2)_r}{(q^2,q^3/b;q^2)_r}\Bigl(\frac{q^2}{b}\Bigr)^r,
\quad f_{2r+1}=0 \quad\text{and}\quad 
\beta_n=\frac{(q^2/b;q^2)_n}{(q,q^2/b;q)_n(q^2;q^2)_n},
\end{equation*}
and
\begin{equation*}
f_{3r}=\frac{(q;q^3)_r}{(q^3;q^3)_r q^r},
\quad f_{3r+1}=f_{3r+2}=0
\quad\text{and}\quad 
\beta_n=\frac{(q;q^3)_n}{(q;q)_n(q;q)_{2n}}
\end{equation*}
both satisfy \eqref{beta}.
By Proposition~\ref{prop1} we therefore have
\begin{multline*}
\sum_{n=0}^{\infty}\frac{(q,q^2/b;q^2)_n q^n}
{(a;q)_{n+1}(q,q/a,q^2/b;q)_n}
-\sum_{r=0}^{\infty}\frac{(q,b;q^2)_r}{(q^2,q^3/b;q^2)_r}
\Bigl(\frac{a^2 q^2}{b}\Bigr)^r \\
=\sum_{n=1}^{\infty} \sum_{r=0}^{\infty}
\frac{(-1)^{n+1} a^n q^{\binom{n-2r}{2}}}{(q,a,q/a;q)_{\infty}}
\frac{(q,b;q^2)_r}{(q^2,q^3/b;q^2)_r}\Bigl(\frac{q^2}{b}\Bigr)^r
(1-q^{(4r+1)n})
\end{multline*}
and
\begin{multline*}
\sum_{n=0}^{\infty}\frac{(q;q^3)_n q^n}{(a;q)_{n+1}(q,q/a;q)_n}
-\frac{(a^3 q^2;q^3)_{\infty}}{(a^3 q;q^3)_{\infty}} \\
=\sum_{n=1}^{\infty}\sum_{r=0}^{\infty}
\frac{(-1)^{n+r+1} a^n q^{r+\binom{n-3r}{2}}}{(q,a,q/a;q)_{\infty}}
\frac{(q;q^3)_r}{(q^3;q^3)_r}
(1-q^{(6r+1)n}),
\end{multline*}
where in the last equation the $q$-binomial theorem \cite[(II.3)]{GR90}
${_1}\phi_0(a;\text{--};q;z)=(az;q)_{\infty}/(z;q)_{\infty}$ has been used
to simplify the second term on the left.

Taking $b=1$ in the quadratic transformation leads to a formula that
might well have been in the lost notebook,
\begin{equation}\label{nNB}
\sum_{n=0}^{\infty}\frac{(q^{n+1};q)_n q^n}
{(a;q)_{n+1}(q^2,q/a;q)_n}=
1+\frac{a+(1+a)\sum_{n=1}^{\infty}(-1)^n a^n q^{\binom{n+1}{2}}}
{(q,a,q/a;q)_{\infty}}.
\end{equation}
For $a=-1$ this further simplifies to the elegant summation
\begin{equation*}
\sum_{n=0}^{\infty}\frac{(q;q^2)_n q^n}{(-q,q^2;q)_n}=
2-\frac{(q;q^2)_{\infty}}{(q^2;q^2)_{\infty}}.
\end{equation*}

Before we continue to derive all of Ramanujan's identities of
the introduction we will slightly change viewpoint and 
reformulate Proposition~\ref{prop1} as a Bailey pair identity.

\section{Partial theta functions and the Bailey lemma}\label{sec4}
Let $\alpha=\{\alpha_n\}_{n=0}^{\infty}$
and $\beta=\{\beta_n\}_{n=0}^{\infty}$. Then the pair of sequences 
$(\alpha,\beta)$ is called a Bailey pair relative to $a$ if \cite{Bailey49}
\begin{equation}\label{BP}
\beta_n=\sum_{r=0}^n \frac{\alpha_r}{(q;q)_{n-r}(aq;q)_{n+r}}.
\end{equation}
Comparing this definition with \eqref{eqprop} and identifying
\begin{equation}\label{alphaf}
\alpha_n=(-1)^n q^{\binom{n}{2}}f_n(1-q^{2n+1})/(1-q)
\end{equation}
we get the following result.
\begin{corollary}\label{cor1}
For $(\alpha,\beta)$ a Bailey pair relative to $q$ there holds
\begin{multline}\label{eqcor}
\sum_{n=0}^{\infty}\frac{\beta_n(q;q)_{2n}q^n}{(a;q)_{n+1}(q/a;q)_n}
-(1-q)\sum_{n=0}^{\infty}
\frac{\alpha_n(-1)^n a^n q^{-\binom{n}{2}}}{1-q^{2n+1}} \\
=\frac{1}{(q^2,a,q/a;q)_{\infty}}\sum_{r=1}^{\infty}
(-1)^{r+1} a^r q^{\binom{r}{2}}
\sum_{n=0}^{\infty}\alpha_n q^{(1-r)n} 
\frac{1-q^{r(2n+1)}}{1-q^{2n+1}},
\end{multline}
provided all sums converge.
\end{corollary}
In view of this result we need to find suitable Bailey pairs
relative to $q$. Before we present many such pairs, we recall
a special case of Bailey's lemma \cite{Andrews84,Paule85}
which states that if $(\alpha,\beta)$ is a Bailey pair relative to $a$,
then so is the pair $(\alpha',\beta')$ given by
\begin{equation*}
\alpha_n'=a^n q^{n^2} \alpha_n \quad\text{and}\quad\beta_n'=
\sum_{r=0}^n \frac{a^r q^{r^2}\beta_r}{(q;q)_{n-r}}.
\end{equation*}
Iterating this leads to what is called the Bailey chain.
As will be shown shortly, all of our theorems of the introduction arise
as Bailey chain identities.
To generate new Bailey pairs we will also use the notion 
of a dual Bailey pair \cite{Andrews81}.
Let $(\alpha,\beta)=(\alpha(a,q),\beta(a,q))$ be a Bailey pair relative
to $a$. Then the pair $(\alpha',\beta')$ given by
\begin{equation*}
\alpha'_n=a^n q^{n^2}\alpha_n(a^{-1},q^{-1}) \quad\text{and}\quad
\beta'_n=a^{-n} q^{-n(n+1)}\beta_n(a^{-1},q^{-1})
\end{equation*}
is again a Bailey pair relative to $a$.
Since in the remainder of this section all Bailey pairs will have $a=q$
we will subsequently drop the phrase ``relative to $q$''.

As our first example we prove Theorem~\ref{T1}.
The required initial Bailey pair is due to Rogers \cite{Rogers17} and given as
item B(3) in Slater's extensive list of Bailey pairs \cite{Slater51},
\begin{equation*}
\alpha_n=(-1)^n q^{n(3n+1)/2}(1-q^{2n+1})/(1-q) \quad\text{and}\quad 
\beta_n=\frac{1}{(q;q)_n}.
\end{equation*}
Iterating this along the Bailey chain one finds
\begin{align*}
\alpha_n^{(k)}&=(-1)^n q^{kn(n+1)+\binom{n}{2}}(1-q^{2n+1})/(1-q) \\
\beta_n^{(k)}&=\sum_{n_1,\dots,n_{k-1}=0}^{\infty}
\frac{q^{N_1^2+\cdots+N_{k-1}^2+N_1+\cdots+N_{k-1}}}
{(q;q)_{n-N_1}(q;q)_{n_1}\cdots (q;q)_{n_{k-1}}}
\end{align*}
for $k$ a positive integer and $(\alpha^{(1)},\beta^{(1)})=
(\alpha,\beta)$. Combining this with Corollary~\ref{cor1} and
applying the triple product identity \eqref{Jtpi} yields
\begin{align*}
\sum_{n=0}^{\infty}&
\frac{\beta_n^{(k)}(q;q)_{2n} q^n}{(a;q)_{n+1}(q/a;q)_n}
-\sum_{n=0}^{\infty} a^n q^{kn(n+1)} \\
&=\frac{1}{(q,a,q/a;q)_{\infty}}\sum_{r=1}^{\infty}
(-1)^{r+1} a^r q^{\binom{r}{2}}
\sum_{n=-\infty}^{\infty}(-1)^n q^{kn(n+1)+\binom{n+1}{2}-rn}  \\
&=\sum_{r=1}^{\infty}(-1)^{r+1} a^r q^{\binom{r}{2}}
\frac{(q^r,q^{2k-r+1},q^{2k+1};q^{2k+1})_{\infty}}{(q,a,q/a;q)_{\infty}}.
\end{align*}
The summand on the right vanishes when $r\equiv 0 \pmod{2k+1}$. Replacing
$r$ by $i+(2k+1)n$ with $i\in\{1,\dots,2k\}$ and $n$ a nonnegative
integer and using \eqref{pred}, we arrive at Theorem~\ref{T1}.

The proof of Theorem \ref{T2} proceeds along the same lines.
We begin with the Bailey pair \cite[E(3)]{Slater51}
\begin{equation*}
\alpha_n=(-1)^n q^{n^2}(1-q^{2n+1})/(1-q) \quad\text{and}\quad
\beta_n=\frac{1}{(q^2;q^2)_n}
\end{equation*}
which implies the iterated pair
\begin{align*}
\alpha_n^{(k)}&=(-1)^n q^{n(kn+k-1)}(1-q^{2n+1})/(1-q) \\
\beta_n^{(k)}&=\sum_{n_1,\dots,n_{k-1}=0}^{\infty}
\frac{q^{N_1^2+\cdots+N_{k-1}^2+N_1+\cdots+N_{k-1}}}
{(q;q)_{n-N_1}(q;q)_{n_1}\cdots(q;q)_{n_{k-2}}
(q^2;q^2)_{n_{k-1}}}.
\end{align*}
Hence, by Corollary~\ref{cor1} and the triple product identity \eqref{Jtpi}, 
\begin{multline*}
\sum_{n=0}^{\infty}\frac{\beta_n^{(k)}(q;q)_{2n}q^n}{(a;q)_{n+1}(q/a;q)_n} \\
=\sum_{n=0}^{\infty}a^n q^{(2k-1)\binom{n+1}{2}}
+\sum_{r=1}^{\infty}(-1)^{r+1} a^r q^{\binom{r}{2}}
\frac{(q^r,q^{2k-r},q^{2k};q^{2k})_{\infty}}{(q,a,q/a;q)_{\infty}}.
\end{multline*}
Replacing $r$ by $i+2kn$ with $i\in\{1,\dots,2k-1\}$ and applying 
\eqref{pred} gives Theorem~\ref{T2}.

Next we turn to Theorem~\ref{T3}, and for the first time a Bailey pair 
not in Slater's list is needed,
\begin{equation}\label{Gnew}
\alpha_n=(-1)^n q^{(3n-1)n/4}(1-q^{2n+1})/(1-q) \quad\text{and}\quad
\beta_n=\frac{1}{(q^2;q^2)_n(-q^{1/2};q)_n}
\end{equation}
which follows from the polynomial identity
\begin{equation}\label{G}
\sum_{j=-\infty}^{\infty}(-1)^j q^{j(3j-1)/4}\qbins{2n+1}{n-j}_q
=(1-q^{2n+1})(q^{1/2};q)_n
\end{equation}
proved in the appendix. We note that by \eqref{alphaf} the
above choice for $\alpha_n$ is equivalent to taking 
$f_r=q^{r(r+1)/4}$ in Proposition~\ref{prop1}.
Substituting the iterated Bailey pair
\begin{align*}
\alpha_n^{(k)}&=(-1)^n q^{(3n-1)n/4+(k-1)n(n+1)}(1-q^{2n+1})/(1-q) \\
\beta_n^{(k)}&=\sum_{n_1,\dots,n_{k-1}=0}^{\infty}
\frac{q^{N_1^2+\cdots+N_{k-1}^2+N_1+\cdots+N_{k-1}}}
{(q;q)_{n-N_1}(q;q)_{n_1}\cdots(q;q)_{n_{k-1}}(-q^{1/2};q^{1/2})_{2n_{k-1}}}
\end{align*}
in equation \eqref{eqcor} gives
\begin{multline*}
\sum_{n=0}^{\infty}\frac{\beta_n^{(k)}(q;q)_{2n}q^n}{(a;q)_{n+1}(q/a;q)_n} 
=\sum_{n=0}^{\infty}a^n q^{(4k-3)(n+1)n/4} \\
+\sum_{r=1}^{\infty}(-1)^{r+1} a^r q^{\binom{r}{2}}
\frac{(q^r,q^{2k-r-1/2},q^{2k-1/2};q^{2k-1/2})_{\infty}}{(q,a,q/a;q)_{\infty}}.
\end{multline*}
Now using that for $g_r$ such that $g_r=0$ if $r\equiv 0\pmod{4k-1}$ there 
holds
\begin{equation*}
\sum_{r=1}^{\infty}g_r=\sum_{i=1}^{2k-1}\sum_{n=0}^{\infty}
\bigl\{g_{(4k-1)n+i}+g_{(4k-1)(n+1)-i}\bigr\}
\end{equation*}
and applying \eqref{pred} we obtain Theorem~\ref{T3}.
For $k=1$ this corresponds to \eqref{R3} with $q\to q^{1/2}$ since
$(q^i,q^{2k-i-1/2},q^{2k-1/2};q^{2k-1/2})_{\infty}=
(q^{1/2};q^{1/2})_{\infty}$ for $k=i=1$.

Finally we prove Theorem~\ref{T4} which is based on the new Bailey pair 
\begin{align}\label{BPA9}
\alpha_n&=(-1)^{\lfloor(4n+1)/3\rfloor}q^{(2n-1)n/3}
\frac{1-q^{2n+1}}{1-q}\chi(n\not\equiv 1\pmod{3}) \\
\beta_n&=\frac{1}{(q;q)_{2n}} \notag
\end{align}
as can be extracted from the polynomial identity
\begin{equation}\label{A9}
\sum_{j=-\infty}^{\infty}\Bigl\{q^{j(6j-1)}\qbins{2n+1}{n-3j}_q-
q^{(2j+1)(3j+1)}\qbins{2n+1}{n-3j}_q\Bigr\}=1-q^{2n+1}
\end{equation}
again proved in the appendix.
Iteration along the Bailey chain yields 
\begin{align*}
\alpha_n^{(k)}&=(-1)^{\lfloor(4n+1)/3\rfloor}q^{(2n-1)n/3+(k-1)n(n+1)}
\frac{1-q^{2n+1}}{1-q}\chi(n\not\equiv 1\pmod{3}) \\
\beta_n^{(k)}&=\sum_{n_1,\dots,n_{k-1}=0}^{\infty}
\frac{q^{N_1^2+\cdots+N_{k-1}^2+N_1+\cdots+N_{k-1}}}
{(q;q)_{n-N_1}(q;q)_{n_1}\cdots(q;q)_{n_{k-2}}(q;q)_{2n_{k-1}}}
\end{align*}
which by Corollary~\ref{cor1} implies
\begin{multline*}
\sum_{n=0}^{\infty}\frac{\beta_n^{(k)}(q;q)_{2n}q^n}{(a;q)_{n+1}(q/a;q)_n}
-\sum_{n=0}^{\infty}
(-1)^n a^{3n} q^{n(3n+1)(2\kappa-3)/2}
(1-a^2 q^{(2\kappa-3)(2n+1)}) \\
=\sum_{r=1}^{\infty}\frac{(-1)^{r+1} a^r q^{\binom{r}{2}}}
{(q,a,q/a;q)_{\infty}}\sum_{n=-\infty}^{\infty} 
q^{\kappa n(3n+1)}\Bigl\{q^{-3nr}-q^{r(3n+1)}\Bigr\},
\end{multline*}
where $\kappa=3k-1$.
Interestingly it is now the quintuple product identity 
\cite[Exercise 5.6]{GR90}
\begin{equation*}
\sum_{n=-\infty}^{\infty}(z^{3n}-z^{-3n-1})q^{n(3n+1)/2}
=(zq,1/z,q;q)_{\infty}(z^2 q,q/z^2;q^2)_{\infty}
\end{equation*}
that is needed to transform the sum over $n$ into a product.
We thus find that the above right-hand side equals
\begin{multline}\label{prod5}
\sum_{r=1}^{\infty}\frac{(-1)^{r+1} a^r q^{\binom{r}{2}}}
{(q,a,q/a;q)_{\infty}}
(q^r,q^{2\kappa-r},q^{2\kappa};q^{2\kappa})_{\infty}
(q^{2\kappa-2r},q^{2\kappa+2r};q^{4\kappa})_{\infty}.
\end{multline}
To rewrite this further we use that for $g_r$ such that $g_r=0$ if
$r\equiv 0 \pmod{\kappa}$ there holds
\begin{equation*}
\sum_{r=1}^{\infty} g_r=\sum_{i=1}^{\kappa-1}\Bigl[g_i+
\sum_{n=1}^{\infty}\bigl(g_{2\kappa n+i}+g_{2\kappa n-i}\bigr)\Bigr].
\end{equation*}
Utilizing this and \eqref{pred} expression \eqref{prod5} becomes
\begin{multline*}
\sum_{i=1}^{\kappa-1}(-1)^{i+1} a^i q^{\binom{i}{2}}
\frac{(q^i,q^{2\kappa-i},q^{2\kappa};q^{2\kappa})_{\infty}
(q^{2\kappa-2i},q^{2\kappa+2i};q^{4\kappa})_{\infty}}{(q,a,q/a;q)_{\infty}} \\
\times \Bigl[1-\sum_{n=1}^{\infty}
a^{2\kappa n-2i}q^{(2\kappa-3)(\kappa n-i)n}
\bigl(1-a^{2i}q^{2(2\kappa-3)in}\bigr)\Bigr],
\end{multline*}
concluding the proof of Theorem~\ref{T4}.

In the remainder of this section we will prove several further 
partial theta series identities that do not reduce to identities
of Ramanujan when $k=1$.
In fact, our first example assumes $k\geq 2$ for reasons of convergence.
Calculating the Bailey pair dual to \eqref{BPA9} gives
\begin{align*}
\alpha_n&=(-1)^{\lfloor(4n+1)/3\rfloor}q^{(n-2)n/3}
\frac{1-q^{2n+1}}{1-q}\chi(n\not\equiv 1\pmod{3}) \\
\beta_n&=\frac{q^{n(n-1)}}{(q;q)_{2n}}.
\end{align*}
Iterating this Bailey pair and copying the previous proof 
one readily finds the following companion to Theorem~\ref{T4}.
\begin{theorem}\label{T5}
For $k\geq 2$, $\kappa=3k-2$ and
$N_j=n_j+n_{j+1}+\cdots+n_{k-1}$,
\begin{align*}
\sum_{n=0}^{\infty}&
\frac{(q;q)_{2n} q^n}{(a;q)_{n+1}(q/a;q)_n}
\sum_{n_1,\dots,n_{k-1}=0}^{\infty}
\frac{q^{N_1^2+\cdots+N_{k-2}^2+2N_{k-1}^2+N_1+\cdots+N_{k-2}}}
{(q;q)_{n-N_1}(q;q)_{n_1}\cdots(q;q)_{n_{k-2}}(q;q)_{2n_{k-1}}} \\
&=\sum_{n=0}^{\infty} (-1)^n a^{3n} q^{(2\kappa-3)(3n+1)n/2}
\bigl(1-a^2 q^{(2\kappa-3)(2n+1)}\bigr) \\
& \quad +\sum_{i=1}^{\kappa-1}(-1)^{i+1} a^i q^{\binom{i}{2}}
\frac{(q^i,q^{2\kappa-i},q^{2\kappa};q^{2\kappa})_{\infty}
(q^{2\kappa-2i},q^{2\kappa+2i};q^{4\kappa})_{\infty}}{(q,a,q/a;q)_{\infty}} \\
&\qquad \qquad \qquad \times \Bigl[1-\sum_{n=1}^{\infty}
a^{2\kappa n-2i}q^{(2\kappa-3)(\kappa n-i)n}
\bigl(1-a^{2i}q^{2(2\kappa-3)in}\bigr)\Bigr].
\end{align*}
\end{theorem}

Next we consider the Bailey pair
\begin{equation}\label{bcs}
\alpha_n=q^{\binom{n}{2}}(1-q^{2n+1})/(1-q)\quad\text{and}\quad
\beta_n=\frac{(-1;q)_n}{(q;q)_n(q^2;q^2)_n}
\end{equation}
which by \eqref{alphaf} corresponds to $f_n=(-1)^n$. Hence \eqref{bcs}
follows by taking $b=-c=q$ in \eqref{f} and \eqref{betabc}. 
Inserting the iterated pair in \eqref{eqcor} gives the next theorem.
\begin{theorem}
For $k$ a nonnegative integer and $N_j=n_j+n_{j+1}+\cdots+n_k$,
\begin{multline*}
\sum_{n=0}^{\infty}
\frac{(q;q)_{2n} q^n}{(a;q)_{n+1}(q/a;q)_n}
\sum_{n_1,\dots,n_k=0}^{\infty}
\frac{q^{N_1^2+\cdots+N_k^2+N_1+\cdots+N_k}(-1;q)_{n_k}}
{(q;q)_{n-N_1}(q;q)_{n_1}\cdots(q;q)_{n_k}(q^2;q^2)_{n_k}}\\
=\sum_{n=0}^{\infty} (-1)^n a^n q^{kn(n+1)}
+\sum_{r=1}^{\infty}\Bigl\{\sum_{n=0}^{\infty}-\sum_{n=-\infty}^{-1}\Bigr\}
\frac{(-1)^{r+1}a^r q^{kn(n+1)+\binom{r-n}{2}}}{(q,a,q/a;q)_{\infty}}.
\end{multline*}
\end{theorem}
Note that the sum over $n$ in the second term on the right takes the
form of a false theta series for which there is no product form.
This can be traced back to the fact that for $f_n=(-1)^n$ there holds 
$f_{-n-1}=-f_n$ instead of the usual $f_n=f_{-n-1}$.
(Note that for generic $b$ and $c$ \eqref{f} satisfies
$f_{-n-1}=f_n$ but that this is not necessarily so if
either $b$ or $c$ assumes the ``singular'' value $q$.)

When $k=0$ the right-hand side trivializes and we get
\begin{equation}\label{k0}
1+2\sum_{n=1}^{\infty}\frac{(q^n;q)_n q^n}{(q,aq,q/a;q)_n}=
\frac{1-a}{1+a}\biggl(1-2\sum_{n=1}^{\infty}
\frac{(-1)^n a^n q^{\binom{n}{2}}}{(q,a,q/a;q)_{\infty}}\biggr)
\end{equation}
which simplifies nicely for $a=1$ to
\begin{equation}\label{a1}
1+2\sum_{n=1}^{\infty}\frac{(q^n;q)_n q^n}{(q;q)^3_n}=
\frac{1}{(q;q)^3_{\infty}}\sum_{n=0}^{\infty}(-1)^n q^{\binom{n+1}{2}}.
\end{equation}

Similar to the previous example we take $b=q$ and $c\to\infty$
in \eqref{f} and \eqref{betabc} to obtain
\begin{equation}\label{bcs2}
\alpha_n=q^{n^2}(1-q^{2n+1})/(1-q)\quad\text{and}\quad
\beta_n=\frac{1}{(q;q)^2_n}.
\end{equation}
Carrying out the usual calculations this yields the second-last
theorem of this section.
\begin{theorem}\label{T6}
For $k$ a nonnegative integer and $N_j=n_j+n_{j+1}+\cdots+n_k$,
\begin{multline*}
\sum_{n=0}^{\infty}
\frac{(q;q)_{2n} q^n}{(a;q)_{n+1}(q/a;q)_n}
\sum_{n_1,\dots,n_k=0}^{\infty}
\frac{q^{N_1^2+\cdots+N_k^2+N_1+\cdots+N_k}}
{(q;q)_{n-N_1}(q;q)_{n_1}\cdots(q;q)_{n_{k-1}}(q;q)^2_{n_k}}\\
=\sum_{n=0}^{\infty} (-1)^n a^n q^{(2k+1)\binom{n+1}{2}}
+\sum_{r=1}^{\infty}\Bigl\{\sum_{n=0}^{\infty}-\sum_{n=-\infty}^{-1}\Bigr\}
\frac{(-1)^{r+1}a^r q^{(2k+1)\binom{n+1}{2}+\binom{r-n}{2}}}
{(q,a,q/a;q)_{\infty}}.
\end{multline*}
\end{theorem}
Again a dramatic simplification occurs for the smallest value of $k$. By
\begin{multline*}
\sum_{r=1}^{\infty}\Bigl\{\sum_{n=0}^{\infty}-\sum_{n=-\infty}^{-1}\Bigr\}
(-1)^{r+1}a^r q^{\binom{n+1}{2}+\binom{r-n}{2}}=
\sum_{r=1}^{\infty}\sum_{n=0}^{r-1}
(-1)^{r+1}a^r q^{\binom{n+1}{2}+\binom{r-n}{2}} \\
=\sum_{n=0}^{\infty}\sum_{r=0}^{\infty}
(-1)^{r+n}a^{r+n+1} q^{\binom{n+1}{2}+\binom{r+1}{2}} 
=a\biggl(\:\sum_{n=0}^{\infty} (-1)^n a^n q^{\binom{n+1}{2}}\biggr)^2,
\end{multline*}
the $k=0$ instance of Theorem~\ref{T6} becomes
\begin{equation}\label{eqk0}
\sum_{n=0}^{\infty} \frac{(q^{n+1};q)_n q^n}{(a;q)_{n+1}(q,q/a;q)_n}
=\sum_{n=0}^{\infty}(-1)^n a^n q^{\binom{n+1}{2}}
+a\frac{\displaystyle
\biggl(\:\sum_{n=0}^{\infty}(-1)^n a^n q^{\binom{n+1}{2}}\biggr)^2}
{(q,a,q/a;q)_{\infty}}.
\end{equation}

Our next Bailey pair can be expressed most concisely as
\begin{equation*}
\alpha_{2n}=\alpha_{-2n-1}=(-1)^n q^{n(3n-1)}(1-q^{4n+1})/(1-q) 
\quad\text{and}\quad
\beta_n=\frac{1}{(q;q)_n(q;q^2)_n},
\end{equation*}
where $n$ in $\alpha_{-2n-1}$ is assumed to be negative.
This Bailey pair can be read off from the polynomial identity
\begin{equation}\label{C8}
\sum_{j=-\infty}^{\infty}(-1)^j q^{j(3j-1)}(1-q^{4j+1})
\qbins{2n+1}{n-2j}_q=(1-q^{2n+1})(-q;q)_n
\end{equation}
established in the appendix. The resulting theorem is similar to those
of the introduction but not quite as beautiful since we cannot carry
out the usual reduction of the sum over $r$. (One can still
use the quasi-periodicity of the triple product on the right
to write the sum over $r$ as a finite sum over $i$ and an infinite
sum over $n$ 
such that the new triple product is $n$-independent, but due to lack of
symmetry, the resulting equation lacks the usual elegance.)
\begin{theorem}
For $k$ a nonnegative integer and $N_j=n_j+n_{j+1}+\cdots+n_k$,
\begin{align*}
\sum_{n=0}^{\infty}&
\frac{(q;q)_{2n} q^n}{(a;q)_{n+1}(q/a;q)_n}
\sum_{n_1,\dots,n_k=0}^{\infty}
\frac{q^{N_1^2+\cdots+N_{k-1}^2+N_1+\cdots+N_{k-1}}}
{(q;q)_{n-N_1}(q;q)_{n_1}\cdots(q;q)_{n_{k-1}}(q;q^2)_{n_{k-1}}}\\
&=\sum_{n=0}^{\infty} (-1)^n a^{2n} q^{n^2+2(k-1)(2n+1)n}
\bigl(1-aq^{(2k-1)(2n+1)}\bigr) \\
&\qquad +\Bigl\{\sum_{r=1}^{\infty}a^r-\sum_{r=-\infty}^{-1}a^{-r}\Bigr\}
(-1)^{r+1}q^{\binom{r}{2}}
\frac{(q^{6k-2r-2},q^{2k+2r},q^{8k-2};q^{8k-2})_{\infty}}
{(q,a,q/a;q)_{\infty}}.
\end{align*}
\end{theorem}
As we have come to expect, the $k=1$ case permits a simplification:
\begin{multline}\label{ft}
\sum_{n=0}^{\infty}
\frac{(-q;q)_n q^n}{(a;q)_{n+1}(q/a;q)_n}
=1+(1+a)\sum_{n=1}^{\infty} (-1)^n a^{2n-1} q^{n^2} \\
+\frac{a-(1+a)\sum_{n=1}^{\infty}a^{3n-1}q^{n(3n-1)/2}(1-aq^n)}
{(q;q^2)_{\infty}(a,q/a;q)_{\infty}}.
\end{multline}
This generalizes the not at all deep, but elegant
\begin{equation*}
\sum_{n=0}^{\infty}\frac{q^n}{(-q;q)_n}=2-\frac{1}{(-q;q)_{\infty}}
\end{equation*}
obtained for $a=-1$.

It will by now be overly clear that the list of nice applications 
of Corollary~\ref{cor1} is sheer endless, and without too much effort
one can obtain many more new Bailey pairs relative to $q$ such that
$\alpha_n$ has a desired factor $(1-q^{2n+1})/(1-q)$.
Most obvious would of course be to use the dual Bailey pairs
corresponding to, for example, \eqref{Gnew} and \eqref{bcs2} (note that
\eqref{bcs} is self-dual), which are given by
\begin{gather*}
\alpha_n=(-1)^n q^{(n-3)n/4}(1-q^{2n+1})/(1-q) \quad\text{and}\quad
\beta_n=\frac{(-1)^n q^{n(n-2)/2}}{(q^2;q^2)_n(-q^{1/2};q)_n} \\
\alpha_n=q^{-n}(1-q^{2n+1})/(1-q)\quad\text{and}\quad
\beta_n=\frac{q^{-n}}{(q;q)^2_n}.
\end{gather*}
The more adventurous reader might also further explore the possibility
of transforming known Bailey pairs into new ones. This has been our
main technique (see the appendix) for deriving new Bailey pairs.
As long as one is willing to allow for Bailey pairs somewhat more
complicated than those presented so far, the number of possible pairs 
appears to be limitless. As an example, we have found numerous 
pairs of the type $\beta_n=q^n/(q;q)_{2n}$ and
\begin{align*}
\alpha_{3n}&=q^{2n(3n-1)}(1-q^{6n+1})/(1-q) \\
\alpha_{3n-1}&=-q^{2n(3n-2)+1}(1-q^{6n-1})/(1-q) \\
\alpha_{3n+1}&=-q^{2n(3n+1)}(1-q^{2n+1})(1-q^{6n+3})/(1-q).
\end{align*}
Unlike our earlier Bailey pairs this follows from a polynomial
identity that can be viewed as a linear combination of
alternating sums over $q$-binomial coefficients,
\begin{multline}\label{A10}
\sum_{j=-\infty}^{\infty}q^{2j(3j-1)}
\Bigl\{\qbins{2n+1}{n-3j}_q-\qbins{2n+1}{n-3j+2}_q\Bigr\} \\
-q\sum_{j=-\infty}^{\infty}q^{2j(3j+2)}
\Bigl\{\qbins{2n+1}{n-3j}_q-\qbins{2n+1}{n-3j-1}_q\Bigr\}
=q^n(1-q^{2n+1}).
\end{multline}
Admittedly the identities obtained after iterating pairs such as
the one above are too involved to be of great interest, but
direct substitution in \eqref{eqcor} often leads to results not
much beyond those of Ramanujan. Our present example, for example,
yields (after some tedious but elementary calculations)
\begin{multline*}
\sum_{n=0}^{\infty}\frac{q^{2n}}{(a;q)_{n+1}(q/a;q)_n}
=1+a+(1+a^2)\sum_{n=1}^{\infty}(-1)^n a^{3n-2}q^{n(3n-1)/2}(1+aq^n) \\
+\frac{a^2+(1+a^2)\sum_{n=1}^{\infty}(-1)^n a^{2n}q^{\binom{n+1}{2}}}
{(a,q/a;q)_{\infty}}.
\end{multline*}
This is so close to \eqref{R0} that it is surprising Ramanujan missed it.
A rather curious formula arises when we set $a=-1$,
\begin{equation*}
\sum_{n=0}^{\infty}\frac{q^{2n}}{(-q;q)^2_n}-
\frac{1+2\sum_{n=1}^{\infty}(-1)^n q^{\binom{n+1}{2}}}
{(-q;q)^2_{\infty}}
=4\sum_{n=1}^{\infty}q^{n(3n-1)/2}(1-q^n).
\end{equation*}
We challenge the reader to explain why all nonzero coefficients on the
right are $\pm 4$.

\section{Residual identities}\label{sec5}
In the following two sections we exploit Andrews' observation 
\cite{Andrews84} that by calculating the residue around the pole $a=q^N$ 
in Ramanujan's partial theta function identities and by then invoking 
analyticity to replace $q^N$ by $a$, new, nontrivial identities arise. 
Of course, instead of considering just Ramanujan's identities we will 
apply Andrews' trick to the more general theorems obtained in 
the previous section.
As a first example however, we treat \eqref{bc} in some detail.
Let $a_0=q^N$ with $N$ a nonnegative integer and multiply both sides of
\eqref{bc} by $(a-a_0)$. The resulting identity is then of the form
\begin{equation*}
\sum_{n=0}^{\infty} f_n(a)-\sum_{r=0}^{\infty} g_r(a)=
\sum_{n=1}^{\infty} h_n(a),
\end{equation*}
with $\lim_{a\to a_0} f_n(a)=f_n(a_0)\chi(n\geq N)$ and
$\lim_{a\to a_0} g_r(a)=0$.
We thus infer that
\begin{equation*}
\sum_{n=0}^{\infty} f_{n+N}(a_0)=\sum_{n=1}^{\infty} h_n(a_0).
\end{equation*}
By straightforward calculations this can be put in the form
\begin{multline}\label{bcres}
\frac{(aq^2/b,aq^2/c,q^2/bc;q)_{\infty}}
{(q^2/b,q^2/c,aq^2/bc;q)_{\infty}}
\sum_{n=0}^{\infty}
\frac{(a^2q^{n+1},aq^2/bc;q)_n q^n}{(q,aq,aq^2/b,aq^2/c;q)_n} \\
=\sum_{n=1}^{\infty}
\frac{(-1)^n a^{n-1} q^{\binom{n+1}{2}}}{(q,q,aq;q)_{\infty}}
\sum_{r=-\infty}^{\infty}
\frac{(-1)^r q^{\binom{r}{2}}(b,c;q)_r}{(q^2/b,q^2/c;q)_r}
\Bigl(\frac{q^{n+3}}{bc}\Bigr)^r,
\end{multline}
where we have replaced $a_0$ by $a$. By standard analyticity arguments
we may now assume $a$ to be an indeterminate.
When $b$ and $c$ tend to infinity the sum over $r$ can be carried out
by the triple product identity leading to \cite[Eq. (8.1)]{Andrews81}
\begin{equation}\label{RA}
\sum_{n=0}^{\infty}
\frac{(a^2 q^{n+1};q)_n q^n}{(q,aq;q)_n}
=\frac{1}{(q,aq;q)_{\infty}}
\sum_{n=0}^{\infty}a^{3n}q^{n(3n+2)}(1-aq^{2n+1}).
\end{equation}
Similarly, when $b=-q$ and $c\to\infty$ we obtain 
\begin{equation}\label{RN}
\sum_{n=0}^{\infty}
\frac{(a q^{n+1};q)_n q^n}{(q;q)_n(aq^2;q^2)_n}
=\frac{1}{(q;q)_{\infty}(aq^2;q^2)_{\infty}}
\sum_{n=0}^{\infty}(-1)^n a^n q^{n(n+1)}.
\end{equation}
The last two identities are in fact closely related to other results
from the lost notebook. If we take $a=-1$ in \eqref{RA} and then
apply Heine's fundamental transformation \cite[Eq. (III.1)]{GR90}
\begin{equation}\label{Heine}
{_2}\phi_1(a,b;c;q,z)=
\frac{(b,az;q)_{\infty}}{(c,z;q)_{\infty}}
{_2}\phi_1(c/b,z;az;q,b),
\end{equation}
with $b=-a=q^{1/2}$, $c=0$ and $z=q$ to transform the left side, we obtain
after replacing $q$ by $q^2$ \cite[p. 37]{Ramanujan88}
\begin{equation*}
\sum_{n=0}^{\infty}\frac{q^n}{(-q;q^2)_{n+1}}=
\sum_{n=0}^{\infty}(-1)^n q^{2n(3n+2)}(1+q^{4n+2}).
\end{equation*}
This is the second ``marvelous'' formula given by Andrews in
his introduction to the lost notebook \cite[Eq. (1.2)]{Andrews79}.
In \cite{Andrews79} Andrews proofs this result by different means.

To transform \eqref{RN} into another of Ramanujan's formulas is 
only slightly more involved. First we use
\begin{equation*}
\frac{(aq^{n+1};q)_n}{(aq^2;q^2)_n}=\frac{(aq;q^2)_n}{(aq;q)_n}=
\frac{(-aq;q)_n(aq;q^2)_n}{(a^2q^2;q^2)_n}
\end{equation*}
to write \eqref{RN} as
\begin{equation*}
\sum_{n=0}^{\infty}\frac{(-aq;q)_n(aq;q^2)_n q^n}{(q;q)_n(a^2q^2;q^2)_n}
=\frac{1}{(q;q)_{\infty}(a q^2;q^2)_{\infty}}
\sum_{n=0}^{\infty}(-1)^n a^n q^{n(n+1)}.
\end{equation*}
By the quadratic analogue of Heine's fundamental transformation
\cite[Thm. A$_3$]{Andrews66}
\begin{equation*}
\sum_{n=0}^{\infty}\frac{(a;q^2)_n(b;q)_{2n}z^n}{(q^2;q^2)_n(c;q)_{2n}}=
\frac{(b;q)_{\infty}(az;q^2)_{\infty}}{(c;q)_{\infty}(z;q^2)_{\infty}}
\sum_{n=0}^{\infty}\frac{(c/b;q)_n(z;q^2)_n b^n}{(q;q)_n(az;q^2)_n},
\end{equation*}
with $a=z\to a q$, $b\to q$ and $c\to -a q^2$ this gives the result
\begin{equation*}
\sum_{n=0}^{\infty}\frac{(q;q^2)_n(aq;q^2)_n(aq)^n}{(-aq;q)_{2n+1}}=
\sum_{n=0}^{\infty}(-1)^n a^n  q^{n(n+1)}.
\end{equation*}
For $a=1$ this is a lost notebook identity \cite[p. 13]{Ramanujan88} proved by 
Andrews~\cite[Eq. (6.2)]{Andrews81}.

A further interesting specialization of \eqref{bcres} arises
when we set $a=-1$. Then the double sum on the right can be
simplified to a single sum resulting in a new $_3\phi_2$ transformation.
First we keep $a$ general and use \cite[Lemma 2]{Andrews81b}
\begin{multline*}
\sum_{n=0}^{\infty} \frac{(q/a;q)_{n+1}(a,b,c;q)_n}{(q;q)_{2n+1}}
\Bigl(\frac{q^2}{bc}\Bigr)^n  \\
=\frac{(q^2/b,q^2/c;q)_{\infty}}{(q,q^2/bc;q)_{\infty}}
\sum_{r=-\infty}^{\infty}
\frac{(-1)^r q^{\binom{r}{2}}(b,c;q)_r}{(q^2/b,q^2/c;q)_r}
\Bigl(\frac{aq^2}{bc}\Bigr)^r,
\end{multline*}
with $a=q^{n+1}$ to write \eqref{bcres} as
\begin{multline*}
\frac{(aq^2/b,aq^2/c;q)_{\infty}}{(aq^2/bc;q)_{\infty}}\sum_{n=0}^{\infty}
\frac{(a^2q^{n+1},aq^2/bc;q)_n q^n}{(q,aq,aq^2/b,aq^2/c;q)_n} \\
=\frac{1}{(q,aq;q)_{\infty}}\sum_{r=0}^{\infty} (b,c;q)_r
\Bigl(\frac{aq^2}{bc}\Bigr)^r \sum_{n=0}^{\infty}
\frac{(-1)^n a^n q^{\binom{n+1}{2}}(q^{2r+2};q)_n}{(q;q)_n}.
\end{multline*}
Here we have changed the order of summation on the right and
shifted $n\to n+r$. For $a=-1$ the sum over $n$ yields
$1/(q;q^2)_{r+1}$ by the $b\to\infty$ limit of the Bailey--Daum sum 
\cite[Eq. (II.9)]{GR90}
\begin{equation*}
{_2}\phi_1(a,b;aq/b;q,-q/b)=\frac{(aq,aq^2/b^2;q^2)_{\infty}}
{(q,q^2)_{\infty}(-q/b,aq/b;q)_{\infty}}.
\end{equation*}
Therefore
\begin{multline*}
\frac{(-q^2/b,-q^2/c;q)_{\infty}}{(-q^2/bc;q)_{\infty}}
\sum_{n=0}^{\infty}\frac{(q;q^2)_n(-q^2/bc;q)_n q^n}{(q,-q^2/b,-q^2/c;q)_n}\\
=\frac{1}{(q^2;q^2)_{\infty}}\sum_{r=0}^{\infty}
\frac{(b,c;q)_r}{(q;q^2)_{r+1}}\Bigl(\frac{-q^2}{bc}\Bigr)^r
\end{multline*}
which can be written as the following nonstandard basic hypergeometric 
transformation
\begin{multline*}
\frac{(-q^2/b,-q^2/c;q)_{\infty}}{(-q^2/bc;q)_{\infty}}
{_3\phi_2}\Bigl[\genfrac{}{}{0pt}{}
{q^{1/2},-q^{1/2},-q^2/bc}{-q^2/b,-q^2/c};q,q\Bigr] \\
=\frac{1}{(1-q)(q^2;q^2)_{\infty}}
{_3\phi_2}\Bigl[\genfrac{}{}{0pt}{}
{q,b,c}{q^{3/2},-q^{3/2}};q,\frac{-q^2}{bc}\Bigr].
\end{multline*}

\section{Garret--Ismail--Stanton-type identities}\label{sec6}
Before calculating more residual identities we review
a recent development in the theory of Rogers--Ramanujan identities 
initiated by Garret, Ismail and Stanton \cite{GIS99} and further 
exploited by Andrews, Knopfmacher and Paule \cite{AKP00}, and 
Berkovich and Paule \cite{BP01,BP01b}. As will be shown later,
calculating the residues of Theorems~\ref{T1}--\ref{T4} provides
a surprising connection between partial theta function identities
and Garret--Ismail--Stanton-type identities.

The famous Rogers--Ramanujan identities are given by
\begin{equation}\label{RR}
\sum_{r=0}^{\infty}\frac{q^{r^2}}{(q;q)_r}=
\frac{1}{(q,q^4;q^5)_{\infty}} \quad\text{and}\quad
\sum_{r=0}^{\infty}\frac{q^{r(r+1)}}{(q;q)_r}=
\frac{1}{(q^2,q^3;q^5)_{\infty}}.
\end{equation}
It is easy to see that the left-hand side of the first (second)
Rogers--Ramanujan identity is the generating function
of partitions of $n$ with difference between parts at least $2$
(and no part equal to $1$).
Schur exploited this fact in one of his proofs of the Rogers--Ramanujan 
identities \cite{Schur17} and introduced two sequences of polynomials 
$\{e_n\}_{n=1}^{\infty}$ and $\{d_n\}_{n=1}^{\infty}$ 
where $e_n$ $(d_n)$ is the generating function corresponding to the
left-side of the first (second) Rogers--Ramanujan identity with the
added condition on the partitions that their largest part does not 
exceed $n-1$.
He then went on to show that both $e_n$ and $d_n$ satisfy the recurrence
\begin{equation}\label{xrec}
x_{n+1}=x_n+q^n x_{n-1}.
\end{equation}
Schur's main result was the following closed form expressions for $e_n$
and $d_n$ 
\begin{subequations}
\begin{align}\label{e}
e_n&=\sum_{j=-\infty}^{\infty}(-1)^j q^{j(5j-1)/2}
\qbin{n}{\lfloor{(n-5j+1)/2}\rfloor}_q \\
d_n&=\sum_{j=-\infty}^{\infty}(-1)^j q^{j(5j-3)/2}
\qbin{n}{\lfloor{(n-5j+2)/2}\rfloor}_q. \label{d}
\end{align}
\end{subequations}
To see that this settles the Rogers--Ramanujan identities
observe that by the triple product identity \eqref{Jtpi} $e_n$ and $d_n$ 
tend to the respective right-hand sides of \eqref{RR} in the large $n$ limit.
Alternative representations for $e_n$ and $d_n$,
probably known to Schur, but first explicitly given by MacMahon
\cite[\S 286 and \S 289]{MacMahon16} are
\begin{equation*}
e_n=\sum_{n=0}^{\infty}q^{r^2}\qbin{n-r}{r}_q \quad\text{and}\quad
d_n=\sum_{n=0}^{\infty}q^{r(r+1)}\qbin{n-r-1}{r}_q.
\end{equation*}
The two polynomial analogues of the Rogers--Ramanujan identities 
obtained by equating the different representations for $e_n$ and 
$d_n$ were first given by Andrews in \cite{Andrews70}.
After this introduction we now come to the beautiful discovery of 
Garret, Ismail and Stanton who found that for $m$ a nonnegative integer
\cite[Eq. (3.5)]{GIS99}
\begin{equation}\label{GIS}
\sum_{r=0}^{\infty}\frac{q^{r(r+m)}}{(q;q)_r}=
\frac{(-1)^m q^{-\binom{m}{2}}d_{m-1}}{(q,q^4;q^5)_{\infty}}-
\frac{(-1)^m q^{-\binom{m}{2}}e_{m-1}}{(q^2,q^3;q^5)_{\infty}}.
\end{equation}
Here $e_{-1}=d_0=0$ and $e_0=d_{-1}=1$ consistent with the recurrence
\eqref{xrec}. For $m=0$ and $m=1$ we of course just find the
first and second Rogers--Ramanujan identity.
A polynomial analogue of \eqref{GIS} was found by Andrews \textit{et al.}
\cite[Prop. 1]{AKP00} in the course of proving \eqref{GIS}
via an extended Engel expansion.
Garret \textit{et al.} also found the inverse of \eqref{GIS} given by
\cite[Thm. 3.1]{GIS99}
\begin{equation}\label{GIS2}
\frac{(q^{3-2i},q^{2i+2},q^5;q^5)_{\infty}}{(q;q)_{\infty}}
=\sum_{j=0}^{\lfloor i/2 \rfloor} (-1)^j q^{2j(j-i)+\binom{j+1}{2}}
\qbin{i-j}{j}_q\sum_{r=0}^{\infty}\frac{q^{r(r+i-2j)}}{(q;q)_r}.
\end{equation}

Next we need a major generalization of \eqref{GIS} due to Berkovich
and Paule \cite{BP01}. As a first step they use the
recurrence \eqref{xrec} to obtain the following ``negative $m$
analogue'' of \eqref{GIS} \cite[Eq. (1.10)]{BP01}
\begin{equation*}
\sum_{r=0}^{\infty}\frac{q^{r(r-m)}}{(q;q)_r}=
\frac{e_m(1/q)}{(q,q^4;q^5)_{\infty}}+
\frac{d_m(1/q)}{(q^2,q^3;q^5)_{\infty}}
\end{equation*}
for $m$ a nonnegative integer. Berkovich and Paule then proceed to
generalize this to the Andrews--Gordon identities given by \cite{Andrews74}
\begin{equation*}
\sum_{n_1,\dots,n_{k-1}=0}^{\infty}
\frac{q^{N_1^2+\cdots+N_{k-1}^2+N_i+\cdots+N_{k-1}}}
{(q;q)_{n_1}\cdots(q;q)_{n_{k-1}}}
=\frac{(q^i,q^{\kappa-i},q^{\kappa};q^{\kappa})_{\infty}}
{(q;q)_{\infty}}
\end{equation*}
for $i\in\{1,\dots,k\}$, $\kappa=2k+1$ and $N_j$ defined as usual.
Before we can give their result we need to define generalizations
of the polynomials $e_n$ and $d_n$ as follows
\begin{multline}\label{X}
X^{(p,p')}_{s,b}(n;q)=X^{(p,p')}_{s,b}(n) \\[2mm]
=\sum_{j=-\infty}^{\infty}
\Bigr\{q^{j(p p'j+p'-ps)}\qbins{n}{(n+s-b)/2-p'j}_q-
q^{(pj+1)(p'j+s)}\qbins{n}{(n-s-b)/2-p'j}_q\Bigl\},
\end{multline}
where $p,p',s,b$ and $n$ are integers such that 
$n+s+b$ is even. The following duality relation 
\cite[Eq. (2.3) and (2.9)]{SW99} will be needed later
\begin{equation}\label{Xdual} 
X_{s,1}^{(p,p')}(n;q)=q^{(n-s+1)(n+s-3)/4} X_{s,1}^{(p'-p,p')}(n;1/q).
\end{equation}
Comparing definition \eqref{X} with \eqref{e} and \eqref{d} shows that 
$e_n=X^{(2,5)}_{2,2+\sigma}(n)$
and $d_n=X^{(2,5)}_{1,2+\sigma}(n)$ where $\sigma\in\{0,1\}$ is fixed by
the condition that $n+\sigma$ is even.
We are now prepared to state the generalization of \eqref{GIS}
as found by Berkovich and Paule \cite[Eq. (3.21)]{BP01}
\begin{multline}\label{eqBP}
\sum_{n_1,\dots,n_{k-1}=0}^{\infty}
\frac{q^{N_1^2+\cdots+N_{k-1}^2+N_{i'}+\cdots+N_{k-1}-m N_1}}
{(q;q)_{n_1}\cdots(q;q)_{n_{k-1}}} \\
=\sum_{\substack{i=1 \\i+i'+m \text{ even}}}^{\kappa-1}
\frac{(q^i,q^{\kappa-i},q^{\kappa};q^{\kappa})_{\infty}}
{(q;q)_{\infty}} X^{(2,\kappa)}_{i,i'}(m;1/q)
\end{multline}
for $i'\in\{1,\dots,k\}$ and $m$ a nonnegative integer.
For $k=2$ and $i'=2$ this is \eqref{GIS} but even for $k=2$ and $i'=1$ 
this is new.
A polynomial analogue of the $i'=k$ case of the
above identity can be found in \cite[Eq. (1.31)]{BP01b}.

After this long introduction into Garret--Ismail--Stanton-type generalizations
of identities of the Rogers--Ramanujan-type we return to our partial
theta function identities and calculate the corresponding
residual identities. First we consider Theorem~\ref{T1} which implies
the identity.
\begin{corollary}
For $k\geq 2$ and $\kappa=2k+1$,
\begin{multline}\label{eqcor2}
\sum_{n=0}^{\infty}
\frac{(a^2 q^{n+1};q)_n q^n}{(q;q)_n}
\sum_{n_1,\dots,n_{k-1}=0}^{\infty}
\frac{q^{N_1^2+\cdots+N_{k-1}^2+N_1+\cdots+N_{k-1}}}
{(aq;q)_{n-N_1}(q;q)_{n_1}\cdots(q;q)_{n_{k-1}}} \\
=\sum_{i=1}^{\kappa-1} (-1)^{i+1} a^{i-1} q^{\binom{i}{2}}
\frac{(q^i,q^{\kappa-i},q^{\kappa};q^{\kappa})_{\infty}}
{(q,q,aq;q)_{\infty}}\sum_{n=0}^{\infty} a^{\kappa n}q^{kn(\kappa n+2i)}.
\end{multline}
\end{corollary}
The identity corresponding to $k=1$ is given by \eqref{RA}.

The similarity between \eqref{eqBP} and \eqref{eqcor2} is quite striking
and in the following we will show how equating coefficients of $a^n$ in
the power series expansion of \eqref{eqcor2} leads to \eqref{eqBP} for $i'=1$. 
As the obvious first step we use \eqref{aq}
and \eqref{pred} to rewrite the above identity as
\begin{multline}\label{idk}
\sum_{n=0}^{\infty}\frac{(a^2 q^{n+1};q)_n q^n}{(q;q)_n}
\sum_{n_1,\dots,n_{k-1}=0}^{\infty}\frac{(aq^{n-N_1+1};q)_{\infty}
q^{N_1^2+\cdots+N_{k-1}^2+N_1+\cdots+N_{k-1}}}
{(q;q)_{n_1}\cdots(q;q)_{n_{k-1}}} \\
=\sum_{r=1}^{\infty}(-1)^{r+1} a^{r-1} q^{\binom{r}{2}}
\frac{(q^r,q^{\kappa-r},q^{\kappa};q^{\kappa})_{\infty}}{(q;q)^2_{\infty}}.
\end{multline} 
Now expanding the left side by \eqref{qbt} and \eqref{qE} gives
\begin{multline*}
\text{LHS}\eqref{idk}\\=
\sum_{j,l,n,n_1,\dots,n_{k-1}=0}^{\infty}
\frac{(-1)^{j+l}a^{j+2l} q^{\binom{j+1}{2}+\binom{l+1}{2}+n(j+l+1)-N_1 j+
\sum_{i=1}^{k-1}N_i(N_i+1)}}
{(q;q)_j(q;q)_l(q;q)_{n-l}(q;q)_{n_1}\cdots(q;q)_{n_{k-1}}}.
\end{multline*}
Shifting $n\to n+l$, $j\to j-2l$ and summing over $n$ using \eqref{qe} yields
\begin{multline*}
\text{LHS}\eqref{idk}=\frac{1}{(q;q)_{\infty}}\sum_{j,l=0}^{\infty}
(-1)^{j+l}a^j q^{\binom{j-l+1}{2}+l(l+1)}\qbin{j-l}{l}_q \\
\times \sum_{n_1,\dots,n_{k-1}=0}^{\infty}
\frac{q^{N_1^2+\cdots+N_{k-1}^2+N_1+\cdots+N_{k-1}-(j-2l)N_1}}
{(q;q)_{n_1}\cdots(q;q)_{n_{k-1}}}.
\end{multline*}
Equating coefficients of $a^j$ with the right side of \eqref{idk} thus gives
\begin{multline}\label{rplus}
\sum_{l=0}^{\lfloor r/2 \rfloor}
(-1)^l q^{l(3l-2r+1)/2}\qbin{r-l}{l}_q 
\sum_{n_1,\dots,n_{k-1}=0}^{\infty}
\frac{q^{N_1^2+\cdots+N_{k-1}^2+N_1+\cdots+N_{k-1}-(r-2l)N_1}}
{(q;q)_{n_1}\cdots(q;q)_{n_{k-1}}} \\
=\frac{(q^{r+1},q^{\kappa-r-1},q^{\kappa};q^{\kappa})_{\infty}}
{(q;q)_{\infty}}
\end{multline}
for $r$ a nonnegative integer.
Comparing this result with \eqref{GIS2} it becomes clear we should now 
invert.  By the connection coefficient formula \cite[Eq. (7.2)]{GIS99}
\begin{equation*}
\frac{H_n(x|q)}{(q;q)_n}=
\sum_{j=0}^{\lfloor n/2 \rfloor} \frac{(1-q^{n-2j+1})q^j}
{(q;q)_j(q;q)_{n-j+1}}
\sum_{l=0}^{\lfloor n/2 \rfloor-j}
\frac{(-1)^l p^{\binom{l+1}{2}}(p;p)_{n-2j-l}}
{(p;p)_l} \frac{H_{n-2j-2l}(x|p)}{(p;p)_{n-2j-2l}}
\end{equation*}
for the $q$-Hermite polynomials and by the $q$-Hermite orthogonality
\cite[Exercise 7.22]{GR90} it follows that 
\begin{equation*}
\frac{h_n}{(q;q)_n}=
\sum_{j=0}^{\lfloor n/2 \rfloor} \frac{(1-q^{n-2j+1})q^j}
{(q;q)_j(q;q)_{n-j+1}}
\sum_{l=0}^{\lfloor n/2 \rfloor-j}
\frac{(-1)^l q^{\binom{l+1}{2}}(q;q)_{n-2j-l}}
{(q;q)_l} \frac{h_{n-2j-2l}}{(q;q)_{n-2j-2l}}
\end{equation*}
for an arbitrary sequence $\{h_n\}_{n=0}^{\infty}$. 
Choosing $h_n=q^{n^2/4}f_n$ this implies the following inversion
\begin{subequations}
\begin{align}\label{gr}
g_r&=\sum_{l=0}^{\lfloor r/2 \rfloor}
(-1)^l q^{l(3l-2r+1)/2}\qbin{r-l}{l}_q f_{r-2l} \\
f_m&=(q;q)_m
\sum_{j=0}^{\lfloor m/2 \rfloor} \frac{(1-q^{m-2j+1})q^{j(j-m+1)}}
{(q;q)_j(q;q)_{m-j+1}}\, g_{m-2j}.
\label{fm}
\end{align}
\end{subequations}
This may also be derived without resorting to $q$-Hermite polynomials
using the $q$-Dougall sum \eqref{qDougall} with $bc=aq=q^{-n}$. 
Since \eqref{rplus} is of the form \eqref{gr} we may rewrite it using
\eqref{fm} to find
\begin{align}\label{inv}
&\sum_{n_1,\dots,n_{k-1}=0}^{\infty}
\frac{q^{N_1^2+\cdots+N_{k-1}^2+N_1+\cdots+N_{k-1}-m N_1}}
{(q;q)_{n_1}\cdots(q;q)_{n_{k-1}}} \\
&=(q;q)_m
\sum_{j=0}^{\lfloor m/2 \rfloor} \frac{(1-q^{m-2j+1})q^{j(j-m+1)}}
{(q;q)_j(q;q)_{m-j+1}}
\frac{(q^{m-2j+1},q^{\kappa+2j-m-1},q^{\kappa};q^{\kappa})_{\infty}}
{(q;q)_{\infty}} \notag \\
&=(q;q)_m
\sum_{\substack{j=1 \\ m+j \text{ odd}}}^{m+1} 
\frac{(1-q^j)q^{(j-m-1)(j+m-3)/4}}
{(q;q)_{(m-j+1)/2}(q;q)_{(m+j+1)/2}}
\frac{(q^j,q^{\kappa-j},q^{\kappa};q^{\kappa})_{\infty}}
{(q;q)_{\infty}}.\notag
\end{align}
The left side coincides with the Berkovich--Paule result \eqref{eqBP}
for $i'=1$. To also show that the above right side agrees with the
right side of \eqref{eqBP} requires some manipulations.
If we denote the summand on the right by $S_j$ then a little calculation
shows that $S_j=S_{-j}$. Since also $S_j=0$ if $j\equiv 0\pmod{\kappa}$ 
or $j>m+1$ we may therefore write
\begin{equation*}
\sum_{\substack{j=1 \\ m+j \text{ odd}}}^{m+1} S_j
=\sum_{\substack{i=1 \\ m+i\text{ odd}}}^{\kappa-1}\sum_{r=-\infty}^{\infty}
S_{2\kappa r+i}.
\end{equation*}
Using this as well as \eqref{pred} we arrive at
\begin{multline*}
\text{LHS}\eqref{inv}=(q;q)_m
\sum_{\substack{i=1 \\ m+i\text{ odd}}}^{\kappa-1}
q^{(i-m-1)(i+m-3)/4}
\frac{(q^i,q^{\kappa-i},q^{\kappa};q^{\kappa})_{\infty}}{(q;q)_{\infty}} \\
\times \sum_{r=-\infty}^{\infty}
\frac{(1-q^{2\kappa r+i})q^{((\kappa-2)(\kappa r+i)-\kappa)r}}
{(q;q)_{(m-2\kappa r-i+1)/2}(q;q)_{(m+2\kappa r+i+1)/2}}.
\end{multline*}
By the easily verified
\begin{equation*}
\frac{(1-q^j)(q;q)_m}{(q;q)_{(m-j+1)/2}(q;q)_{(m+j+1)/2}}=
\qbins{m}{(m-j+1)/2}_q-q^j\qbins{m}{(m-j-1)/2}_q
\end{equation*}
it thus follows that
\begin{multline*}
\text{LHS}\eqref{inv}=
\sum_{\substack{i=1 \\ m+i\text{ odd}}}^{\kappa-1} q^{(i-m-1)(i+m-3)/4}
\frac{(q^i,q^{\kappa-i},q^{\kappa};q^{\kappa})_{\infty}}{(q;q)_{\infty}} \\
\times \sum_{r=-\infty}^{\infty} 
q^{((\kappa-2)(\kappa r+i)-\kappa)r}
\Bigl\{\qbins{m}{(m-2\kappa r-i+1)/2}_q-q^{2\kappa r+i}
\qbins{m}{(m-2\kappa r-i-1)/2}_q\Bigr\}.
\end{multline*}
In the first term within the curly braces we replace $r\to -r$ and 
use $\qbins{m+n}{m}=\qbins{m+n}{n}$. 
Comparing the resulting expression with \eqref{X} we find that 
the second line of the above equation is 
$X^{(\kappa-2,\kappa)}_{i,1}(m)$.
By the duality relation \eqref{Xdual} we therefore find
\begin{equation*}
\text{LHS}\eqref{inv}=\sum_{\substack{i=1\\ m+i\text{ odd}}}^{\kappa-1} 
\frac{(q^i,q^{\kappa-i},q^{\kappa};q^{\kappa})_{\infty}}{(q;q)_{\infty}}
X^{(2,\kappa)}_{i,1}(m;1/q)
\end{equation*}
in agreement with the right side of \eqref{eqBP} for $i'=1$.

Before we continue to treat the other theorems of sections~\ref{sec1} and 
\ref{sec4} let us briefly comment on \eqref{rplus}. Because of the
ocurrence of the $q$-binomial coefficient it is clear 
that for negative values of $r$ the left side trivially vanishes. 
The right side,
on the other hand, is nonvanishing for any integer $r$ as long as
$r+1\not\equiv 0 \pmod{\kappa}$, and is in fact symmetric under the
transformation $r\to \kappa-r-2$.
To obtain an identity valid for all integers $r$ one can use the
symmetry of the right side (or more precisely, the quasi-periodicity
under the transformation $r\to -r-2$) to prove that
\begin{multline*}
\sum_{l=-\infty}^{\lfloor r/2 \rfloor}
(-1)^l q^{l(3l-2r+1)/2}\qbin{r-l}{r-2l}_q
\sum_{n_1,\dots,n_{k-1}=0}^{\infty}
\frac{q^{N_1^2+\cdots+N_{k-1}^2+N_1+\cdots+N_{k-1}-(r-2l)N_1}}
{(q;q)_{n_1}\cdots(q;q)_{n_{k-1}}} \\
=\frac{(q^{r+1},q^{\kappa-r-1},q^{\kappa};q^{\kappa})_{\infty}}
{(q;q)_{\infty}},
\end{multline*}
where $r$ is now an arbitrary integer and where the $q$-binomial coefficient
is redefined as $\qbins{m+n}{m}=(q^{n+1};q)_m/(q;q)_m$ for $m$ a nonnegative 
integer and zero otherwise. Note that this implies that the lower bound 
in the sum over $l$ may be optimized to $\min\{0,r+1\}$.

The other theorems on partial theta functions may be applied
in a similar manner. Since each time the calculations only marginally differ,
we will leave out the details and give the most important equations only.
The residual equation corresponding to Theorem~\ref{T2} can be stated
as follows.
\begin{corollary}
For $k\geq 2$ and $\kappa=2k$,
\begin{multline*}
\sum_{n=0}^{\infty}
\frac{(a^2 q^{n+1};q)_n q^n}{(q;q)_n}
\sum_{n_1,\dots,n_{k-1}=0}^{\infty}
\frac{q^{N_1^2+\cdots+N_{k-1}^2+N_1+\cdots+N_{k-1}}}
{(aq;q)_{n-N_1}(q;q)_{n_1}\cdots(q;q)_{n_{k-2}}(q^2;q^2)_{n_{k-1}}} \\
=\sum_{i=1}^{\kappa-1} (-1)^{i+1} a^{i-1} q^{\binom{i}{2}}
\frac{(q^i,q^{\kappa-i},q^{\kappa};q^{\kappa})_{\infty}}
{(q,q,aq;q)_{\infty}}
\sum_{n=0}^{\infty} (-1)^n a^{\kappa n}q^{(\kappa-1)(kn+i)n}.
\end{multline*}
\end{corollary}
The identity corresponding to $k=1$ is given by \eqref{RN}.

By equating coefficients of $a^n$ one finds
\begin{multline*}
\sum_{l=0}^{\lfloor r/2 \rfloor}
(-1)^l q^{l(3l-2r+1)/2}\qbin{r-l}{l}_q
\sum_{n_1,\dots,n_{k-1}=0}^{\infty}
\frac{q^{N_1^2+\cdots+N_{k-1}^2+N_1+\cdots+N_{k-1}-(r-2l)N_1}}
{(q;q)_{n_1}\cdots(q;q)_{n_{k-2}}(q^2;q^2)_{n_{k-1}}} \\
=\frac{(q^{r+1},q^{\kappa-r-1},q^{\kappa};q^{\kappa})_{\infty}}
{(q;q)_{\infty}}
\end{multline*}
for $r$ a nonnegative integer. Inversion of this equation gives
the $i'=1$ instance of
\begin{multline*}
\sum_{n_1,\dots,n_{k-1}=0}^{\infty}
\frac{q^{N_1^2+\cdots+N_{k-1}^2+N_{i'}+\cdots+N_{k-1}-m N_1}}
{(q;q)_{n_1}\cdots(q;q)_{n_{k-2}}(q^2;q^2)_{n_{k-1}}} \\
=\sum_{\substack{i=1\\ m+i+i'\text{ odd}}}^{\kappa-1} 
\frac{(q^i,q^{\kappa-i},q^{\kappa};q^{\kappa})_{\infty}}{(q;q)_{\infty}}
X^{(2,\kappa)}_{i,i'}(m;1/q),
\end{multline*}
where we remind the reader that $\kappa=2k$.
In view of \eqref{eqBP} it is not difficult to guess that the above is true
for all $i'\in\{1,\dots,k\}$.

The Garret--Ismail--Stanton-type identity associated to Theorem~\ref{T3}
takes a slightly different form. First we calculate the residual identity.
\begin{corollary}\label{cor4}
For $k\geq 2$ and $\kappa=2k-1/2$, 
\begin{align*}
\sum_{n=0}^{\infty} & \frac{(a^2 q^{n+1};q)_n q^n}{(q;q)_n}
\sum_{n_1,\dots,n_{k-1}=0}^{\infty}
\frac{q^{N_1^2+\cdots+N_{k-1}^2+N_1+\cdots+N_{k-1}}}
{(aq;q)_{n-N_1}(q;q)_{n_1}\cdots(q;q)_{n_{k-1}}(-q^{1/2};q^{1/2})_{2n_{k-1}}}\\
&\quad =\sum_{i=1}^{2k-1} (-1)^{i+1} a^{i-1} q^{\binom{i}{2}}
\frac{(q^i,q^{\kappa-i},q^{\kappa};q^{\kappa})_{\infty}}
{(q,q,aq;q)_{\infty}} \\
& \qquad \qquad \qquad \times
\sum_{n=0}^{\infty} (-1)^n a^{2\kappa n}q^{2(\kappa-1)(\kappa n+i)n}
\bigl(1+a^{2\kappa-2i}q^{2(\kappa-1)(\kappa-i)(2n+1)}\bigr).
\end{align*}
\end{corollary}
The equation corresponding to $k=1$ is
\begin{equation*}
\sum_{n=0}^{\infty}\frac{(a;q)_n q^n}{(q,a^2;q)_n}
=\frac{(a;q)_{\infty}}{(q,a^2;q)_{\infty}}
\sum_{n=0}^{\infty} (-1)^n a^{3n}q^{(3n-1)n/2}(1+aq^n),
\end{equation*}
where we have replaced $a$ by $a/q^{1/2}$.
By \eqref{Heine} with $b\to 0$, $c\to a^2$ and $z\to q$ this can be 
transformed into
\begin{equation*}
\sum_{n=0}^{\infty}\frac{(-1)^n a^{2n}q^{\binom{n}{2}}}{(a;q)_{n+1}}
=\sum_{n=0}^{\infty} (-1)^n a^{3n}q^{(3n-1)n/2}(1+aq^n)
\end{equation*}
which for $a=q^{1/2}$ yields Rogers' false theta function identity
\cite[p. 333; Eq. (4)]{Rogers17}
\begin{equation*}
\sum_{n=0}^{\infty}\frac{(-1)^n q^{n(n+1)}}{(q;q^2)_{n+1}}=
1+\sum_{n=1}^{\infty} (-1)^n q^{3n^2}(q^{2n}-q^{-2n}).
\end{equation*}
Returning to the more general case, we equate powers of $a^n$ in
Corollary~\ref{cor4} to obtain
\begin{align*}
\sum_{l=0}^{\lfloor r/2 \rfloor} (-1)^l q^{l(3l-2r+1)/2}\qbin{r-l}{l}_q
\sum_{n_1,\dots,n_{k-1}=0}^{\infty}
\frac{q^{N_1^2+\cdots+N_{k-1}^2+N_1+\cdots+N_{k-1}-(r-2l)N_1}}
{(q;q)_{n_1}\cdots(q;q)_{n_{k-1}}(-q^{1/2};q^{1/2})_{2n_{k-1}}}\\
=\frac{(q^{r+1},q^{\kappa-r-1},q^{\kappa};q^{\kappa})_{\infty}}
{(q;q)_{\infty}}.
\end{align*}
After inversion this gives the $i'=1$ instance of
\begin{multline*}
\sum_{n_1,\dots,n_{k-1}=0}^{\infty}
\frac{q^{N_1^2+\cdots+N_{k-1}^2+N_{i'}+\cdots+N_{k-1}-m N_1}}
{(q;q)_{n_1}\cdots(q;q)_{n_{k-1}}(-q^{1/2};q^{1/2})_{2n_{k-1}}}\\
=\sum_{\substack{i=1\\ m+i+i'\text{ odd}}}^{2\kappa-1}
\frac{(q^i,q^{\kappa-i},q^{\kappa};q^{\kappa})_{\infty}}{(q;q)_{\infty}}
X^{(4,2\kappa)}_{i,i'}(m;1/q).
\end{multline*}
Again we conjecture this to hold for all $i'\in\{1,\dots,k\}$ 
and $\kappa=2k-1/2$.
Note the subtle difference with earlier cases in that
the sum over $i$ on the right exceeds $\kappa$.
Since $X_{i,i'}^{(4,2\kappa)}(m;1/q)$ (for $i\in\{1,\dots 2\kappa-1\}$ and 
$i'\in\{1,\dots,k\}$) is a polynomial with only
nonnegative coefficients \cite{ABBBFV87} this implies that for 
$i>\kappa$ the summand on the right (as a power series in $q$) has
nonpositive coefficients.

Finally we treat Theorems~\ref{T4} and \ref{T5} together.
\begin{corollary}\label{C3}
For $\sigma\in\{0,1\}$, $k\geq 2$ and $\kappa=3k-\sigma-1$,
\begin{multline*}
\sum_{n=0}^{\infty}
\frac{(a^2 q^{n+1};q)_n q^n}{(q;q)_n}
\sum_{n_1,\dots,n_{k-1}=0}^{\infty}
\frac{q^{N_1^2+\cdots+N_{k-1}^2+N_1+\cdots+N_{k-1}+\sigma N_{k-1}(N_{k-1}-1)}}
{(aq;q)_{n-N_1}(q;q)_{n_1}\cdots(q;q)_{n_{k-2}}(q;q)_{2n_{k-1}}} \\
=\sum_{i=1}^{\kappa-1}(-1)^{i+1} a^{i-1}q^{\binom{i}{2}}
\frac{(q^i,q^{2\kappa-i},q^{2\kappa};q^{2\kappa})_{\infty}
(q^{2\kappa-2i},q^{2\kappa+2i};q^{4\kappa})_{\infty}}{(q,q,aq;q)_{\infty}} \\
\times \Bigl[1-\sum_{n=1}^{\infty}
a^{2\kappa n-2i}q^{(2\kappa-3)(\kappa n-i)n}
\Bigl\{1-a^{2i}q^{2(2\kappa-3)in}\Bigr\}\Bigr].
\end{multline*}
\end{corollary}
The identity corresponding to $\sigma=0$ and $k=1$ turns out to be a 
special case of Heine's transformation \eqref{Heine} and has therefore 
been omitted. Equating the coefficients of $a^n$ yields
\begin{multline*}
\sum_{l=0}^{\lfloor r/2 \rfloor}(-1)^l q^{l(3l-2r+1)/2}\qbin{r-l}{l}_q \\
\times \sum_{n_1,\dots,n_{k-1}=0}^{\infty}
\frac{q^{N_1^2+\cdots+N_{k-1}^2+N_1+\cdots+N_{k-1}+
\sigma (N^2_{k-1}-N_{k-1})-(r-2l)N_1}}
{(q;q)_{n_1}\cdots(q;q)_{n_{k-2}}(q;q)_{2n_{k-1}}} \\
=\frac{(q^{r+1},q^{2\kappa-r-1},q^{2\kappa};q^{2\kappa})_{\infty}
(q^{2\kappa-2r-2},q^{2\kappa+2r+2};q^{4\kappa})_{\infty}}{(q;q)_{\infty}},
\end{multline*}
provided $r$ is a nonnegative integer. The inversion of this equation gives
the $i'=1$ case of
\begin{multline*}
\sum_{n_1,\dots,n_{k-1}=0}^{\infty}
\frac{q^{N_1^2+\cdots+N_{k-1}^2+N_{i'}+\cdots+N_{k-1}+
\sigma N_{k-1}(N_{k-1}-1)-m N_1}}
{(q;q)_{n_1}\cdots(q;q)_{n_{k-2}}(q;q)_{2n_{k-1}}} \\
=\sum_{\substack{i=1\\ m+i+i'\text{ odd}}}^{\kappa-1} 
\frac{(q^i,q^{2\kappa-i},q^{2\kappa};q^{2\kappa})_{\infty}
(q^{2\kappa-2i},q^{2\kappa+2i},q^{4\kappa};q^{4\kappa})_{\infty}}
{(q;q)_{\infty}} X^{(3,\kappa)}_{i,i'}(m;1/q).
\end{multline*}
Yet again we conjecture this to be true for all $i'\in\{1,\dots,k\}$.

Since the outcomes are less spectacular, we leave the calculation of the
residual identities of the remaining theorems of section~\ref{sec4} to the
reader. Instead we just list some of the simplest cases which we hope are
of some interest.

Calculating the residue around $a=q^N$ in \eqref{nNB} results in
\begin{equation*}
\sum_{n=0}^{\infty}\frac{(a^2 q^{n+1};q)_n q^n}{(q,aq,aq^2;q)_n}=
\frac{1+(1+1/a)\sum_{n=1}^{\infty}(-1)^n a^n q^{\binom{n+1}{2}}}
{(1-q)(q,aq,aq^2;q)_{\infty}}
\end{equation*}
which for $a=1$ simplifies to
\begin{equation*}
\sum_{n=0}^{\infty}\frac{(q^{n+1};q)_n q^n}{(q,q,q^2;q)_n}=
\frac{1+2\sum_{n=1}^{\infty}(-1)^n q^{\binom{n+1}{2}}}
{(q;q)^3_{\infty}}
\end{equation*}
reminiscent of \eqref{a1}.

The residual identity corresponding to \eqref{k0} is
\begin{equation*}
\sum_{n=0}^{\infty}\frac{(a^2;q)_{2n} q^n}{(q,a,aq,a^2q;q)_n}=
\frac{1}{(q,aq,aq;q)_{\infty}}
\sum_{n=0}^{\infty}(-1)^n a^n q^{\binom{n+1}{2}}
\end{equation*}
which generalizes \eqref{a1} obtained for $a=1$. 
This last result becomes more interesting if we compare it with the
analogous result obtained from \eqref{eqk0}:
\begin{equation*}
\sum_{n=0}^{\infty} \frac{(a^2 q^{n+1};q)_n q^n}{(q,aq,aq;q)_n}
=\frac{1}{(q,aq,aq;q)_{\infty}}
\biggl(\:\sum_{n=0}^{\infty}(-1)^n a^n q^{\binom{n+1}{2}}\biggr)^2.
\end{equation*}
Noting the similarity of the above two right-hand sides and using basic
hypergeometric notation we infer that
\begin{equation*}
\biggl({_3\phi_2}\Bigl[\genfrac{}{}{0pt}{}
{-a,aq^{1/2},-aq^{1/2}}{aq,a^2 q};q,q\Bigr]\biggr)^2
=\frac{1}{(q,aq,aq;q)_{\infty}}
{_3\phi_2}\Bigl[\genfrac{}{}{0pt}{}
{-aq,aq^{1/2},-aq^{1/2}}{aq,a^2 q};q,q\Bigr].
\end{equation*}

Finally, calculating the residue around $a=q^N$ in \eqref{ft} yields
\begin{equation*}
\sum_{n=0}^{\infty}
\frac{(-aq;q)_n q^n}{(q,a^2 q;q)_n}
=\frac{(-aq;q)_{\infty}}{(q,a^2q;q)_{\infty}}
\Bigl[1-(1+a)\sum_{n=1}^{\infty}a^{3n-2}q^{n(3n-1)/2}(1-aq^n)\Bigr].
\end{equation*}
By Heine's transformation \eqref{Heine} with $a\to -aq$, $b\to 0$,
$c\to a^2 q$ and $z\to q$ this becomes
\begin{equation*}
\sum_{n=0}^{\infty}\frac{(-1)^n a^{2n} q^{\binom{n+1}{2}}}{(-aq;q)_{n+1}}
=1-(1+a)\sum_{n=1}^{\infty}a^{3n-2}q^{n(3n-1)/2}(1-aq^n).
\end{equation*}
For $a=1$ this yields 
\begin{equation}\label{fth}
\sum_{n=0}^{\infty}\frac{(-1)^n q^{n(n+1)}}{(-q^2;q^2)_{n+1}}
=1+2\sum_{n=1}^{\infty}q^{3n^2}(q^n-q^{-n}).
\end{equation}
Although this is a false theta function identity 
not in Rogers' paper it readily follows that 
$\eqref{fth}-\eqref{f6}=\eqref{f6}-1$ with \cite[p. 333; Eq. (6)]{Rogers17}
\begin{equation}\label{f6}
\sum_{n=0}^{\infty}\frac{(-1)^n q^{n(n+1)}}{(-q^2;q^2)_n}
=1+\sum_{n=1}^{\infty}q^{3n^2}(q^n-q^{-n}).
\end{equation}

\section{Discussion}\label{sec7}
In the abstract we stated that many of Ramanujan's partial theta function
identities can be generalized by the method developed in this paper.
In the main text, however, we restricted ourselves to Ramanujan's
identities \eqref{R0} and \eqref{R1}--\eqref{R3}, which all have a very
similar structure dictated by Proposition~\ref{prop1}. 
To conclude we will give one example of how simple
modifications lead to generalizations of other partial
theta function identities of the lost notebook.
The identity we will generalize in our example is \cite[p. 37]{Ramanujan88}
\begin{multline*}
\sum_{n=0}^{\infty}\frac{q^{2n+1}}{(aq,q/a;q^2)_{n+1}} \\
=\sum_{n=0}^{\infty}(-1)^{n+1}a^{3n+1} q^{n(3n+2)}(1+aq^{2n+1})
+\sum_{n=0}^{\infty}\frac{(-1)^n a^{2n+1}q^{n(n+1)}}
{(aq,q/a;q^2)_{\infty}}
\end{multline*}
which was first proved by Andrews \cite[Eq. (3.9)]{Andrews81}.

First we take $b=q^2/a$ in \eqref{ab} and divide both sides by
$(1-q/a)$. This yields
\begin{multline}\label{R4}
\sum_{n=0}^{\infty}\frac{(q;q)_{2n+1}q^n}{(a,q/a;q)_{n+1}}
\sum_{r=0}^n \frac{(-1)^r q^{\binom{r}{2}}f_r(1-q^{2r+2})}
{(q;q)_{n-r}(q;q)_{n+r+2}}
+\sum_{r=0}^{\infty}(a/q)^{r+1} f_r \\
=\frac{1}{(q,a,q/a;q)_{\infty}}\sum_{n=1}^{\infty}\Bigl\{
\sum_{r=0}^{\infty}f_r+q^{-n}\sum_{r=-\infty}^{-1}f_{-r-1}\Bigr\}
(-1)^{n+r} a^n q^{\binom{n+r+1}{2}}
\end{multline}
provided all sums converge.
Recalling definition \eqref{BP} of a Bailey pair this can be restated 
as follows.
\begin{corollary}\label{cor2b}
For $(\alpha,\beta)$ a Bailey pair relative to $q^2$ there holds
\begin{multline}\label{eqcor2b}
\sum_{n=0}^{\infty}\frac{\beta_n (q^2;q)_{2n}q^n}{(a,q/a;q)_{n+1}}
+(1-q^2)\sum_{n=0}^{\infty}\frac{\alpha_n (-1)^n (a/q)^{n+1}
q^{-\binom{n}{2}}} {1-q^{2n+2}} \\
=\frac{(1-q^2)}{(q,a,q/a;q)_{\infty}}\sum_{r=1}^{\infty}
(-1)^{r+1} (a/q)^r q^{\binom{r}{2}}
\sum_{n=0}^{\infty}\alpha_n q^{(1-r)n}\frac{1-q^{r(2n+2)}}{1-q^{2n+2}}.
\end{multline}
\end{corollary}
The Bailey pair required to turn this into \eqref{R4} is
\begin{align}\label{BPA11}
\alpha_n&=(-1)^{\lfloor 4n/3\rfloor}q^{n(2n+1)/3}
\frac{1-q^{2n+2}}{1-q^2} \chi(n\not\equiv 2\pmod{3})\\
\beta_n&=\frac{1}{(q^2;q)_{2n}}. \notag
\end{align}
The proof of this pair comes down to the proof of the polynomial identity
\begin{equation}\label{A11}
\sum_{j=-\infty}^{\infty}\Bigl\{q^{j(6j+1)}\qbins{2n+2}{n-3j}_q-
q^{(2j+1)(3j+1)}\qbins{2n+2}{n-3j-1}_q\Bigr\}=1-q^{2n+2}
\end{equation}
established in the appendix.
If we insert \eqref{BPA11} in \eqref{eqcor2b}, replace
$a$ by $aq^{1/2}$ followed by $q\to q^2$ we find \eqref{R4}.
If first we use the Bailey lemma to obtain the iterated Bailey pair
\begin{align*}
\alpha_n&=(-1)^{\lfloor 4n/3\rfloor}q^{n(2n+1)/3+(k-1)n(n+2)}
\frac{1-q^{2n+2}}{1-q^2} 
\chi(n\not\equiv 2\pmod{3})\\
\beta_n&=\sum_{n_1,\dots,n_{k-1}=0}^{\infty}
\frac{q^{N_1^2+\cdots+N_{k-1}^2+2N_1+\cdots+2N_{k-1}}}
{(q;q)_{n-N_1}(q;q)_{n_1}\cdots (q;q)_{n_{k-2}}(q^2;q)_{2n_{k-1}}}
\end{align*}
relative to $q^2$ and insert this in \eqref{eqcor2b} with $a$
replaced by $a q^{1/2}$, we obtain our final theorem. 
\begin{theorem}
For $k\geq 1$, $\kappa=3k-1$ and $N_j=n_j+n_{j+1}+\cdots+n_{k-1}$, 
\begin{align*}
\sum_{n=0}^{\infty}&\frac{(q^2;q)_{2n}q^{n+1/2}}{(aq^{1/2},q^{1/2}/a;q)_{n+1}}
\sum_{n_1,\dots,n_{k-1}=0}^{\infty}
\frac{q^{N_1^2+\cdots+N_{k-1}^2+2N_1+\cdots+2N_{k-1}}}
{(q;q)_{n-N_1}(q;q)_{n_1}\cdots (q;q)_{n_{k-2}}(q^2;q)_{2n_{k-1}}} \\
&=\sum_{n=0}^{\infty}(-1)^{n+1} a^{3n+1}q^{(2\kappa-3)(3n+2)n/2}
\bigl(1+a q^{(2\kappa-3)(2n+1)/2}\bigr) \\
&\quad +\sum_{i=1}^{\kappa-1}(-1)^{i+1} a^i q^{(i-1)^2/2}
\frac{(q^{\kappa+i},q^{\kappa-i},q^{2\kappa};q^{2\kappa})_{\infty}
(q^{2i},q^{4\kappa-2i};q^{4\kappa})_{\infty}}
{(q,aq^{1/2},q^{1/2}/a;q)_{\infty}} \\
&\qquad \qquad \qquad \times \sum_{n=0}^{\infty}a^{2\kappa n}
q^{(2\kappa-3)(\kappa n+i)n}
\Bigl\{1-a^{2\kappa-2i}q^{(2\kappa-3)(\kappa-i)(2n+1)}\Bigr\}.
\end{align*}
\end{theorem}

As a final comment we should note that it does not seem possible to
prove and generalize all of Ramanujan's partial theta function formulas
using \eqref{ab}, and it is an open problem whether one can modify
our approach to extend an identity like \cite[p. 39]{Ramanujan88}
\begin{multline*} 
\sum_{n=0}^{\infty} \frac{q^{3n^2}}{(a;q^3)_{n+1}(q^3/a;q^3)_n}- 
q\sum_{n=1}^{\infty}\frac{q^{3n(n-1)}}{(aq,q^2/a;q^3)_n}+
q/a\sum_{n=1}^{\infty}\frac{q^{3n(n-1)}}{(q/a,aq^2;q^3)_n} \\
=\frac{(q;q)^2_{\infty}}{(q^3;q^3)_{\infty}(a,q/a;q)_{\infty}}
\end{multline*}
which contains partial and complete theta functions of different moduli.

\appendix

\section{Proofs of polynomial identities}
In this appendix we prove the various polynomial identities used in the
main text for extracting Bailey pairs. All proofs are based on identities
obtained by Rogers in his classic 1917 paper \cite{Rogers17} on 
Rogers--Ramanujan-type identities, which we transform using
the $q$-binomial recurrences
\begin{equation*}
\qbins{m+n}{m}_q=
\qbins{m+n-1}{m}_q+q^{n-m}\qbins{m+n-1}{m-1}_q 
=\qbins{m+n-1}{m-1}_q+q^m \qbins{m+n-1}{m}_q. 
\end{equation*}

For the proof of \eqref{G} we require 
\begin{gather}\label{G1}
\sum_{j=-\infty}^{\infty}(-1)^j q^{j(3j-1)/4}\qbins{2n}{n-j}_q
=(q^{1/2};q)_n \\ \label{G3}
\sum_{j=-\infty}^{\infty}(-1)^j q^{3j(j-1)/4}\qbins{2n}{n-j}_q
=q^n(q^{1/2};q)_n
\end{gather}
equivalent to the Bailey pairs G(1) and G(3) in Slater's list \cite{Slater51}.
Applying the first $q$-binomial recurrence to the left-hand side
of \eqref{G} yields
\begin{align*}
\text{LHS}\eqref{G}&=\text{LHS}\eqref{G1}+q^{n+1}
\sum_{j=-\infty}^{\infty}(-1)^j q^{3(j+1)j/4}\qbins{2n}{n-j-1}_q \\
&=(q^{1/2};q)_n-q^{n+1}\text{LHS}\eqref{G3}=(1-q^{2n+1})(q^{1/2};q)_n,
\end{align*}
where the second equality follows after the variable change $j\to j-1$.

The ingredients needed for the proofs of \eqref{A9}, \eqref{A10} and 
\eqref{A11} are polynomial identities equivalent to the
Bailey pairs A(1)--A(4) \cite{Slater51},
\begin{gather}\label{A1}
\sum_{j=-\infty}^{\infty}\Bigl\{
q^{j(6j-1)}\qbins{2n}{n-3j}_q-
q^{(2j+1)(3j+1)}\qbins{2n}{n-3j-1}_q\Bigr\}=1 \\
\label{A2}
\sum_{j=-\infty}^{\infty}\Bigl\{
q^{j(6j+1)}\qbins{2n+1}{n-3j}_q-
q^{(2j+1)(3j+1)}\qbins{2n+1}{n-3j-1}_q\Bigr\}=1 \\
\label{A3}
\sum_{j=-\infty}^{\infty}q^{j(6j+2)}\Bigl\{
\qbins{2n}{n-3j}_q-\qbins{2n}{n-3j-1}_q\Bigr\}=q^n \\
\label{A4}
\sum_{j=-\infty}^{\infty}q^{j(6j+4)}\Bigl\{
\qbins{2n+1}{n-3j}_q-\qbins{2n+1}{n-3j-1}_q\Bigr\}=q^n.
\end{gather}
To show \eqref{A9} is true we take its left-hand side
and apply the first (second) $q$-binomial recurrence to the 
first (second) term of the summand to find
\begin{equation*}
\text{LHS}\eqref{A9}=\text{LHS}\eqref{A1}-q^{n+1}\text{LHS}\eqref{A3}
=1-q^{2n+1}.
\end{equation*}
To establish \eqref{A11} we take its left-hand side and apply the 
first $q$-binomial recurrences to both terms of the summand to find
\begin{align*}
\text{LHS}\eqref{A11}&=\text{LHS}\eqref{A2}-q^{n+2}
\sum_{j=-\infty}^{\infty}\Bigl\{q^{(j+1)(6j+2)}\qbins{2n+1}{n-3j-2}_q-
q^{j(6j+4)}\qbins{2n+1}{n-3j-1}_q\Bigr\} \\
&=1-q^{n+2}\text{LHS}\eqref{A4}=1-q^{2n+2},
\end{align*}
where the second equality follows by the variable change $j\to -j-1$
in the first term of the sum over $j$ and the symmetry 
$\qbins{m+n}{m}_q=\qbins{m+n}{n}_q$.
The polynomial identity \eqref{A10} requires some more work.
First observe that the second line on the left of \eqref{A10} is precisely
$-q\text{LHS}\eqref{A4}=-q^{n+1}$. Using the 
first (second) $q$-binomial recurrence on the first (second) term of the 
summand on the first line (and making some trivial variable changes) thus
gives
\begin{multline*}
\text{LHS}\eqref{A10}=\text{LHS}\eqref{A3}-q^{n+1} \\
+q^{n+1}\sum_{j=-\infty}^{\infty}\Bigr\{q^{j(6j+1)}\qbins{2n}{n-3j-1}_q-
q^{(2j+1)(3j+1)}\qbins{2n}{n-3j-2}_q\Bigr\}.
\end{multline*}
Next we expand the remaining sum over $j$ using the first $q$-binomial
recurrence. This leads to
\begin{align*}
\text{LHS}\eqref{A10}&=q^n-q^{n+1} 
+q^{n+1}\text{LHS}\eqref{A2}|_{n\to n-1} \\
& \qquad -q^{2n+2}\sum_{j=-\infty}^{\infty}
\Bigr\{q^{(j+1)(6j+2)}\qbins{2n-1}{n-3j-3}_q
-q^{j(6j+4)}\qbins{2n-1}{n-3j-2}_q\Bigr\} \\
&=q^n-q^{2n+2}\text{LHS}\eqref{A4}|_{n\to n-1}=
q^n(1-q^{2n+1}).
\end{align*}

Finally we deal with \eqref{C8} for which we need three polynomial
identities equivalent to the Bailey pairs 
C(1), C(2) and C(4) \cite{Slater51},
\begin{gather}\label{C1}
\sum_{j=-\infty}^{\infty}(-1)^j q^{j(3j-1)}\qbins{2n}{n-2j}_q=(-q;q)_n \\
\label{C2}
\sum_{j=-\infty}^{\infty}(-1)^j q^{j(3j+1)}
\Bigl\{\qbins{2n}{n-2j}_q-\qbins{2n}{n-2j-1}_q\Bigr\}=q^n(-q;q)_n \\
\label{C4}
\sum_{j=-\infty}^{\infty}(-1)^j q^{3j(j+1)}
\qbins{2n+1}{n-2j}_q=q^n(-q;q)_n.
\end{gather}
First note that the second sum in \eqref{C8} (corresponding to the term 
$q^{4j+1}$ in $(1-q^{4j+1})$) is $-q\,\text{LHS}\eqref{C4}=-q^{n+1}(-q;q)_n$.
We therefore need to show that the first sum, which will be denoted
by $S_1$, equals $(1+q^{n+1}(1-q^n))(-q;q)_n$.
Now, by the first $q$-binomial recurrence, 
\begin{align*}
S_1&=\sum_{j=-\infty}^{\infty}(-1)^j q^{j(3j-1)}\qbins{2n}{n-2j}_q+
q^{n+1}\sum_{j=-\infty}^{\infty}(-1)^j q^{j(3j+1)}\qbins{2n}{n-2j-1}_q \\
&=\text{LHS}\eqref{C1}+q^{n+1}(\text{LHS}\eqref{C1}-\text{LHS}\eqref{C2})=
(1+q^{n+1}(1-q^n))(-q;q)_n.
\end{align*}
Here we note that in calculating $\text{LHS}\eqref{C1}-\text{LHS}\eqref{C2}$
one should first replace $j\to -j$ in \eqref{C1} and use 
$\qbins{m+n}{m}=\qbins{m+n}{n}$.

\bibliographystyle{amsplain}

\end{document}